\documentclass[12pt,a4paper,reqno]{amsart}   
\usepackage{ifthen} 
\newboolean{colorinverted}
\setboolean{colorinverted}{true} 
\setboolean{colorinverted}{false} 
 
\ifthenelse{\boolean{colorinverted}}{ 
\usepackage[usenames,dvipsnames,svgnames,table]{xcolor} 
\pagecolor{black} 
\color{Gray} 
}

\usepackage{geometry}
\usepackage{mathtools}   
\usepackage{url}
 
\usepackage[myheadings]{fullpage}  
\usepackage[english]{babel} 
\usepackage{enumerate}   
  
\usepackage[full]{textcomp} 
\usepackage{newtxtext} 
\usepackage{cabin} 
\usepackage{zlmtt}
\usepackage[bigdelims,vvarbb]{newtxmath} 
\usepackage[cal=boondoxo]{mathalfa} 

\usepackage{microtype} %
  
\newcommand{\eps}{\varepsilon}

\newcommand{\LL}{\mathcal{L}}

\newcommand{\Q}{\mathbb{Q}}

\newcommand{\R}{\mathbb{R}} 

\newcommand{\Z}{\mathbb{Z}}

\newcommand{\Odi}[1]{\Odip{}{#1}}

\newcommand{\Odip}[2]{\mathcal{O}_{#1}\left(#2\right)}

\newcommand{\odip}[2]{{o}_{#1}\!\left(#2\right)}
\newcommand{\odi}[1]{\odip{}{#1}}

\newcommand{\ULI}{\text{ULI}}
\newcommand{\LLI}{\text{LLI}} 

\renewcommand{\qedsymbol}{$\square$}
\newenvironment{Proof}[1][Proof]{\par\noindent\textbf{#1.}~}
{\hfill\qedsymbol\smallskip\par}

\newtheorem{Theorem}{Theorem}
\newtheorem{Lemma}{Lemma}

\newtheoremstyle{remark}
{10pt}
{6pt}
{\rm} 
{}
{\bfseries}
{.}
{.5em}
{\thmname{#1}\thmnumber{ #2}\thmnote{ (#3)}}
 \theoremstyle{remark}

\usepackage{caption} 
\newcommand{\bound}{10^7}
\newcommand{\boundtable}{1000}
\newcommand{\indirizzoweb}{\url{https://www.math.unipd.it/~languasc/LS_zero.html}}
\newcommand{\clusteraddress}{\url{https://hpc.math.unipd.it}}

\newcommand{\constonemin}{0.0124862668}
\newcommand{\constoneminattained}{7105733}
\newcommand{\constonemax}{0.6267599041}
\newcommand{\constonemaxattained}{23}
\newcommand{\consttwomin}{0.0091904477}
\newcommand{\consttwominattained}{7105733}
\newcommand{\consttwomax}{0.4206022969}
\newcommand{\consttwomaxattained}{311}

\allowdisplaybreaks 
\begin{document} 

\title[Numerical estimates on the Landau-Siegel zero]
{Numerical estimates on the Landau-Siegel zero \\
and other related quantities}  
 \author{Alessandro Languasco}

\begin{abstract}   
Let $q$ be an odd prime, $\chi$ be a non-principal  Dirichlet character $\bmod\ q$
and $L(s,\chi)$ be the associated Dirichlet $L$-function.
For every odd prime $q\le\bound$,
we show that 
\(
 L(1,\chi_\square)  > c_{1} \log q
\)
and 
\(
\beta < 1- \frac{c_{2}}{\log q},
\)
where $c_1=\constonemin\dotsc$,  $c_2=\consttwomin\dotsc$, $\chi_{\square}$ 
is the quadratic Dirichlet character $\bmod\ q$ 
and $\beta\in (0,1)$ is the Landau-Siegel zero, if it exists,
of the set of such Dirichlet $L$-functions.
As a by-product of the computations here performed,
we also obtained some information about  Littlewood's and Joshi's
bounds on $L(1,\chi_\square)$   and on the class number of 
the imaginary quadratic field $\Q(\sqrt{-q})$.
\end{abstract}
\maketitle
\makeatletter
\def\subsubsection{\@startsection{subsubsection}{3}%
  \z@{.3\linespacing\@plus.5\linespacing}{-.5em}%
  {\normalfont\bfseries}} 
\makeatother
\section{Introduction}
Let $q$ be an odd prime,  $\chi$ be  a non-principal Dirichlet character $\bmod\ q$
and $L(s,\chi)$ be the associated Dirichlet $L$-function.
Exploiting a fast algorithm to compute $\vert L(1,\chi) \vert$  
developed in  \cite{Languasco2021} and \cite{Languasco2023a},
we will  obtain the values of $L(1,\chi_\square)$  
for every quadratic Dirichlet character $\chi_\square$  $\bmod\ q$, and
for every odd prime $q \le \bound$. %
Using these informations, we will derive some numerical estimates on the Landau-Siegel zero, 
if it exists, attached to the set
\begin{equation}
\label{L-set-def}
\LL:=\{ L(s,\chi_\square)\colon \chi_\square  \bmod{q}, 3\le q \le \bound, q\ \textrm{prime}\}.
\end{equation}
In this way we will obtain  the following result.
\begin{Theorem}
\label{Thm-beta}
Let $q$ be an odd prime, $q\le \bound$, let $\chi_\square$ be 
the quadratic Dirichlet character mod $q$ and let
$\beta:=\beta_\LL\in (0,1)$ be the Landau-Siegel zero of $\LL$, if it exists.
We have that there exist two computable constants
$c_{1},c_{2}>0$ such that
\begin{equation}
\label{c1-def}
  L(1,\chi_\square)  > c_{1} \log q
\end{equation}
for every $L(s,\chi_\square) \in \LL$, 
and
\begin{equation}
\label{c2-def}
\beta < 1- \frac{c_{2}}{\log q}.
\end{equation}
 \end{Theorem}
The values of $c_{1}$ and $c_{2}$ are obtained as the minimal ones of the analogous
quantities we computed for each $q$; more precisely, for every odd prime $q \le \bound$  we got 
the values $c_1(q)$ and $c_2(q)$ such that \eqref{c1-def}-\eqref{c2-def}  hold
and then we obtained 
\[
c_{1} = \min_{\substack{3\le q \le \bound \\ q\, \textrm{prime}}} 
c_{1}(q)  = \constonemin\dotsc
\quad
\textrm{and}
\quad
c_{2} =  \min_{\substack{3\le q \le \bound \\ q\, \textrm{prime}}}  
c_{2}(q)  = \consttwomin\dotsc
\]
both attained at $q=\constoneminattained$ ($\chi_\square$ is even).
Analogously we can get the less important quantities
\[
C_{1} = \max_{\substack{3\le q \le \bound \\ q\, \textrm{prime}}}  
c_{1}(q) =  \constonemax\dotsc
\quad
\textrm{and}
\quad
C_{2} =  \max_{\substack{3\le q \le \bound \\ q\, \textrm{prime}}}  
c_{2}(q)  = \consttwomax \dotsc,
\] 
attained at  $q=\constonemaxattained$, 
respectively $q=311$ (in both cases $\chi_\square$ is odd).
The dependence of $c_{2}(q)$ from $c_{1}(q)$ is stated  in \eqref{c1-c2-function} below
and it involves also a computed value for 
$S(q) : = \sum_{n=2}^{q} (\log n)/n$ and an effective version of the 
P\'olya-Vinogradov inequality proved by Lapkova \cite{Lapkova2018} in 2018.

Recalling that  Watkins \cite{Watkins2004} showed that there are no Siegel zeroes for $q\le 3\cdot 10^8$  whenever
$\chi_\square$ is odd, and Platt \cite{Platt2015} reached the same conclusion for $\chi_\square$ even and $q\le 4\cdot 10^5$,
our Theorem \ref{Thm-beta} is meaningful only for $4\cdot 10^5 < q \le \bound$, $q$ prime, and $\chi_\square$ 
even. Moreover, we also recall that Morrill and Trudgian \cite{MorrillT2020} proved that for $q\ge 3$, the 
function $L(s,\chi_\square)$ has at most one real zero $\beta$
with $1-0.933/\log q < \beta <1$.

\begin{figure}[ht] 
\includegraphics[scale=0.5,angle=0]{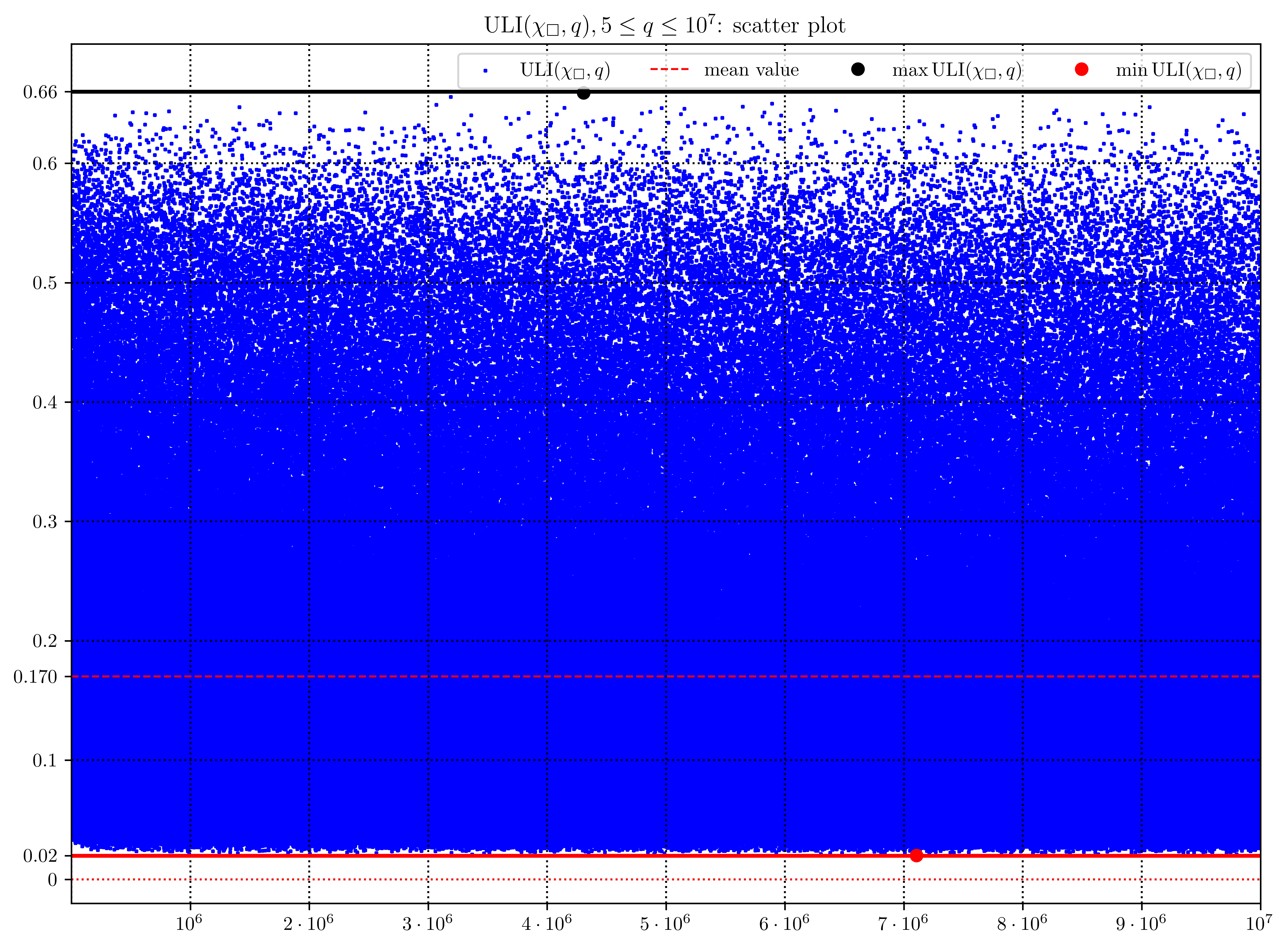}  
\caption{\small{The values of $\ULI(\chi_\square,q)$, $q$ prime, $5\le q\le \bound$.
$\min \ULI(\chi_\square,q)  = 0.0200472032\dotsc$ attained at 
$q=7105733$; 
$ \max \ULI(\chi_\square,q) =  0.6590147671\dotsc$ attained at 
$q=4305479$.
The black line corresponds to $0.66$;  the red one  to $0.02$.
$\ULI(\chi_\square,q)$ is defined in \eqref{ULI-LLI-def}.
The red dashed line corresponds to the mean value.
 }}
\label{fig-ULI} 
 \end{figure}

\begin{figure}[ht] 
\includegraphics[scale=0.5,angle=0]{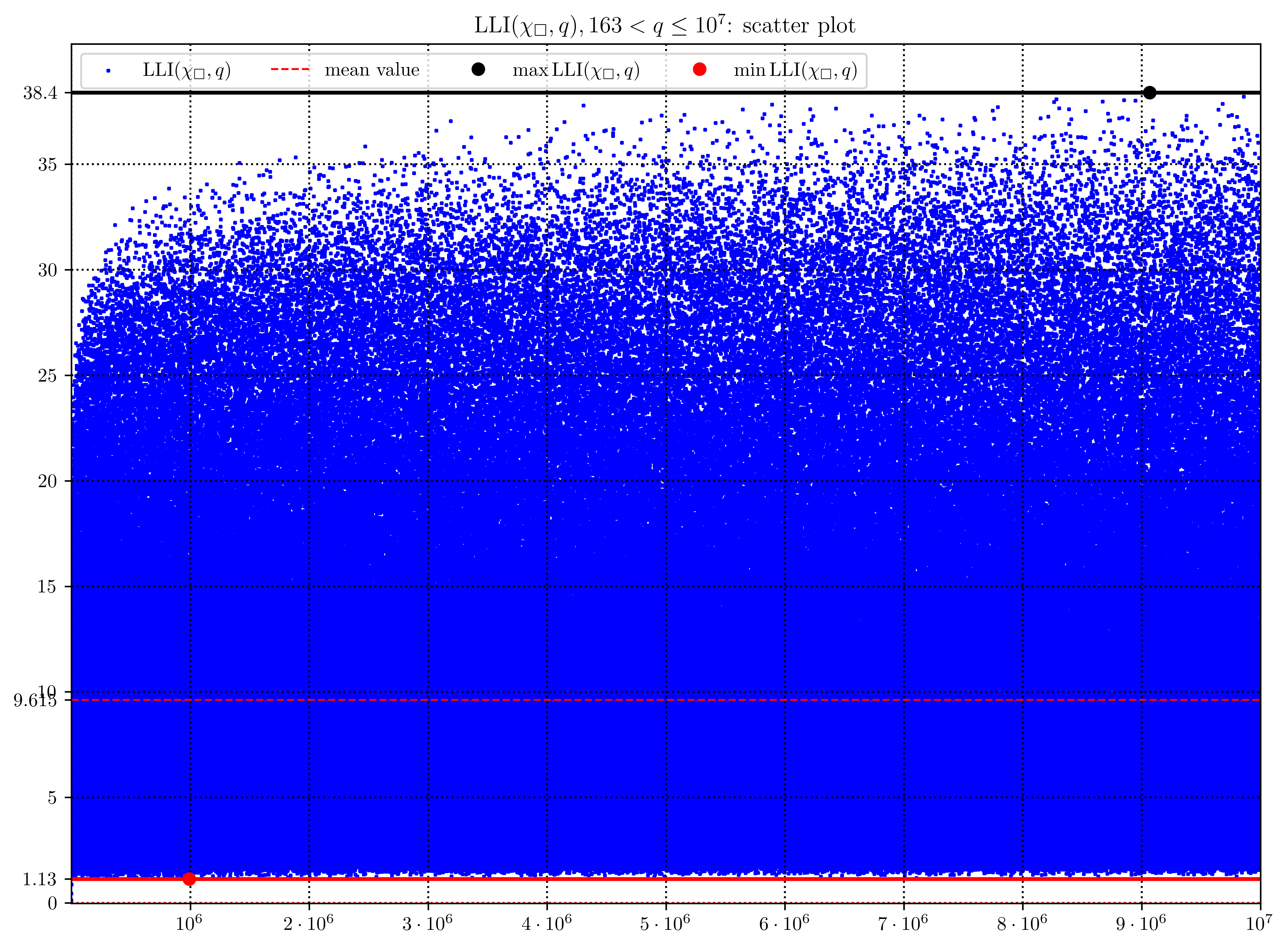}  
\caption{\small{The values of $\LLI(\chi_\square,q)$, $q$ prime, $163< q\le \bound$. 
$\min \LLI(\chi_\square,q)  = 1.1302203128\dotsc$ attained at 
$q=991027$; 
$ \max \LLI(\chi_\square,q) =  38.3973766224\dotsc$ attained at 
$q=9067439$.
The black line corresponds to $38.4$;  the red one  to $1.13$.
$\LLI(\chi_\square,q)$ is defined in \eqref{ULI-LLI-def}.
The red dashed line corresponds to the mean value.
 }}
\label{fig-LLI} 
 \end{figure}

\begin{figure}[ht] 
\includegraphics[scale=0.5,angle=0]{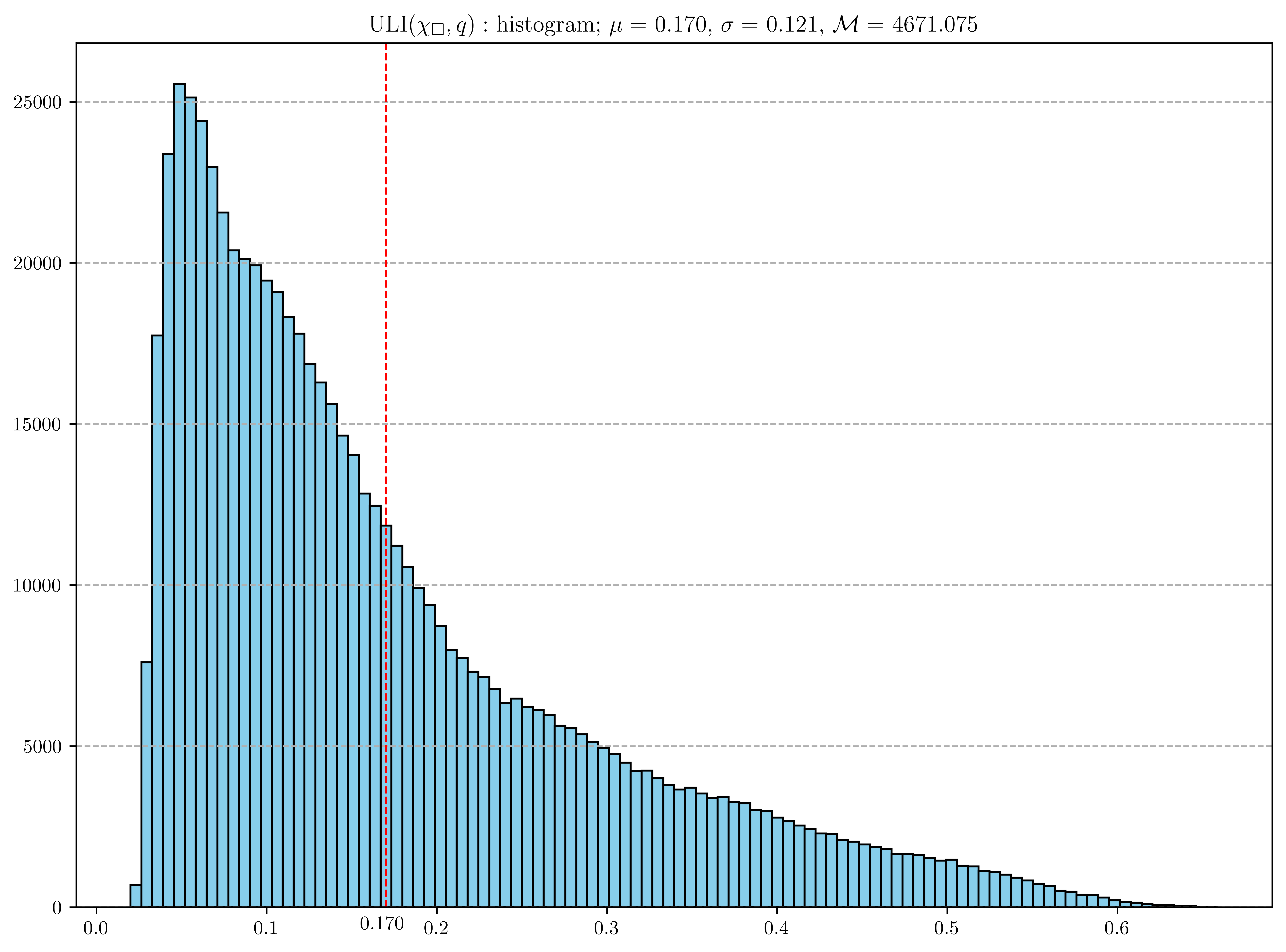}  
\caption{\small{Histogram about the values of $\ULI(\chi_\square,q)$, $q$ prime, $5\le q\le \bound$.
Intervals length $:= I =0.0070286432\dotsc$;
number of primes $5\le q\le \bound := {\mathcal P} = 664577$;
mass $:=  \mathcal{M} =  I \cdot {\mathcal P}$;
mean $:= \mu = 0.1701388804\dotsc$; 
standard deviation $:= \sigma = 0.1213103612\dotsc$.
$\ULI(\chi_\square,q)$ is defined in \eqref{ULI-LLI-def}.
The red dashed line corresponds to the mean value.
 }}
\label{fig-ULI-histo} 
\end{figure}

\begin{figure}[ht] 
\includegraphics[scale=0.5,angle=0]{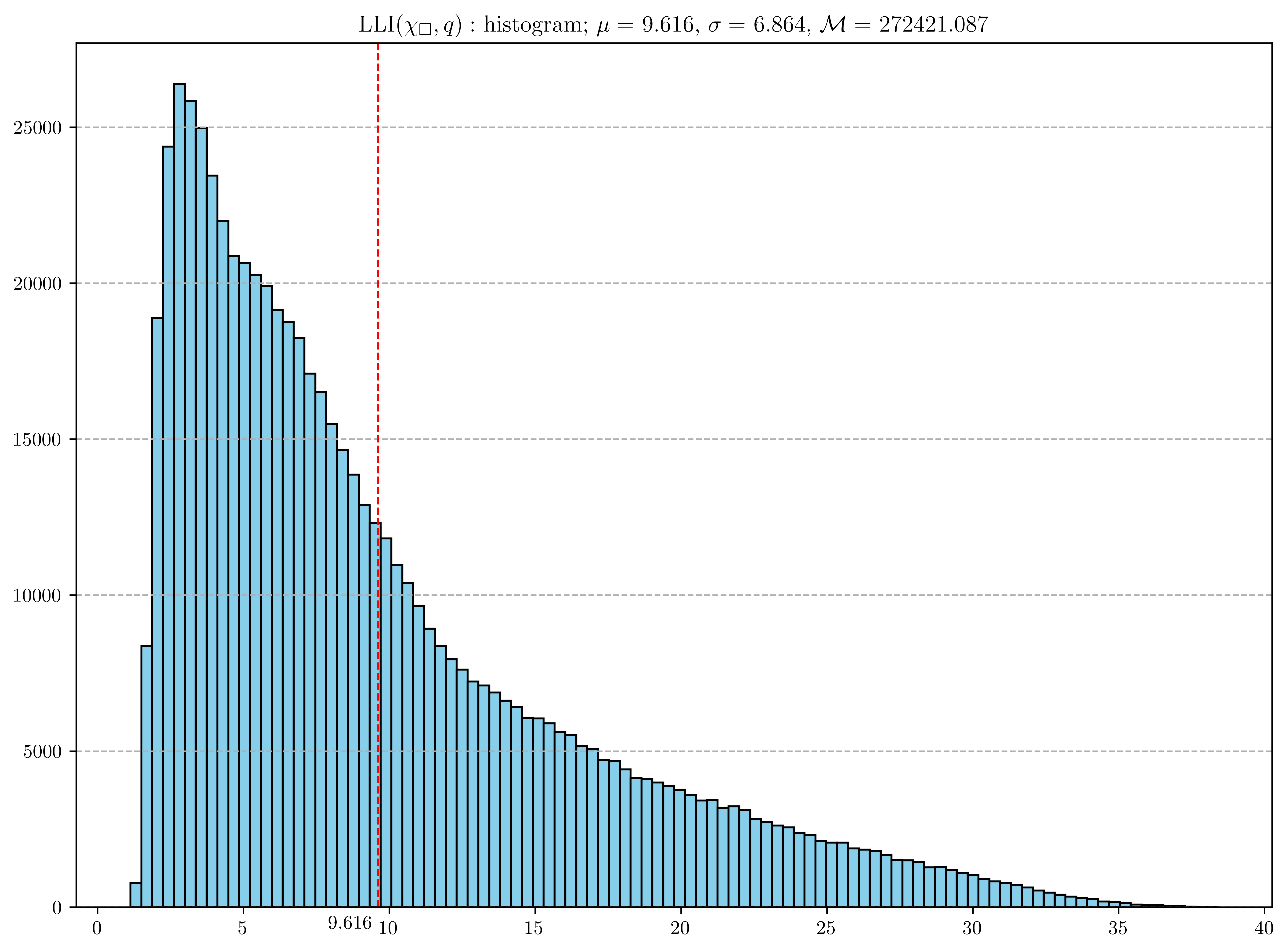}  
\caption{\small{Histogram about the values of $\LLI(\chi_\square,q)$, $q$ prime,  $163 < q\le \bound$.
Intervals length $:= I = 0.4099387194\dotsc$;
number of primes $163 < q\le \bound := {\mathcal P} = 664541$;
mass $:=  \mathcal{M} =  I \cdot {\mathcal P}$;
mean $:= \mu = 9.6157071784\dotsc$; 
standard deviation $:= \sigma = 6.8639742972\dotsc$.
$\LLI(\chi_\square,q)$ is defined in \eqref{ULI-LLI-def}.
The red dashed line corresponds to the mean value.
}}
\label{fig-LLI-histo} 
\end{figure}

A second set of results is about the size of $L(1,\chi_\square)$. First of all, 
we recall that, assuming the Riemann Hypothesis for $L(s,\chi_\square)$,
Littlewood \cite{Littlewood1928} proved in 1928 that
\begin{equation}
\label{Littlewood-bounds}
\Bigl(
\frac{12e^\gamma}{\pi^2}(1+\odi{1})\log \log  q
\Bigr)^{-1}
<
L(1,\chi_\square)
<
2e^\gamma (1+\odi{1}) \log \log  q
\end{equation}
as $q$ tends to infinity, where $\gamma$ is the Euler-Mascheroni constant.
In 1973 Shanks \cite{Shanks1973} numerically studied 
the behaviour of the \emph{upper} and \emph{lower Littlewood indices} defined as
\begin{equation}
\label{ULI-LLI-def}
\ULI (d,\chi_d):=
\frac{ L(1,\chi_d)}{2e^\gamma\log \log \vert d\vert}
\quad
\text{and}
\quad
\LLI(d,\chi_d) := L(1,\chi_d)\frac{12e^\gamma}{\pi^2}\log \log \vert d\vert
\end{equation}
for several small discriminants $d$. 
 Such computations were extended by 
Williams-Broere \cite{WilliamsB1976} in 1976 and by 
Jacobson-Ramachandran-Williams \cite{Jacobson2006} in 2006.
From the values of $L(1,\chi_\square)$ here obtained we can infer that
\begin{equation*}
1.130
<
\LLI(q,\chi_\square) 
<
38.398
\end{equation*}
for every prime $7\le q \le \bound$ (the  minimal 
value  $1.1302203128\dotsc$ is attained 
at $q=991027$, $\chi_\square$ is odd;
the maximal value $38.3973766224\dotsc$ is attained at 
$q=9067439$, $\chi_\square$ is odd) with the unique exception of $q=163$ 
($\LLI(163,\chi_\square) = 0.8675157625\dotsc$, $\chi_\square$ is odd).
The only cases in which $\LLI(q,\chi_\square) <1$ were attained at $q=3,5,163$.
Moreover,  we also have that
\begin{equation*}
0.020
<
\ULI(q,\chi_\square) 
<
0.660
\end{equation*}
for every  prime $5 \le q \le \bound$ (the maximal value 
$0.6590147671\dotsc$ is attained at $q=4305479$, $\chi_\square$ is odd;
the minimal value $0.0200472032\dotsc$ is attained at 
$q=7105733$, $\chi_\square$ is even).
The only case in which $\ULI(q,\chi_\square) >1$ was attained at $q=3$.
In Figures \ref{fig-ULI}-\ref{fig-LLI} we inserted two scatter plots about $\ULI$ and $\LLI$
and in Figures \ref{fig-ULI-histo}-\ref{fig-LLI-histo} we presented the corresponding
histograms.

Moreover, in 1970  Joshi \cite[Theorem 1]{Joshi1970} proved that 
\begin{equation}
\label{Joshi-first}
\limsup_{q\to \infty} \frac{L(1, \chi_\square)}{\log\log q} \geq e^{\gamma} 
\end{equation}
and
\begin{equation}
\label{Joshi-second}
\liminf_{q\to \infty} L(1, \chi_\square) \log\log q \le  \frac{\pi^2}{6e^{\gamma}}.
\end{equation}
Our computation shows that for only 
about  $1.96$\%  of the  odd primes 
$q\le \bound$ we have that $L(1, \chi_\square)/\log\log q  \geq e^{\gamma}$;
to be more precise, such an inequality holds for exactly $13036$ odd primes $q\le \bound$ 
($\pi(10^7)=664579$).
The first ten cases are attained at $q=3, 7, 71, 191, 239, 311, 479, 719, 839, 1151$.
Moreover, for about  $2.17$\%  of the  odd primes 
$q\le \bound$ we have that $L(1, \chi_\square) \log\log q \le  \pi^2/(6 e^{\gamma})$;
to be more precise, such an inequality holds for exactly $14453$ odd primes $q\le \bound$.
The first ten cases are attained at $q=3, 5, 7, 11, 13, 19, 29, 43, 53, 67$.

Let now $q\ge 5$ and  denote as  $h(-q)$ the class number of the imaginary quadratic field $\Q(\sqrt{-q})$.
The famous Dirichlet class number formula, \emph{i.e.},
\begin{equation}
\label{Dirichlet-class-number}
h(-q) = (\sqrt{q}/\pi) \, L(1,\chi_\square),
\end{equation} 
gives us the chance to derive the values of $h(-q)$ from the computed values of $L(1,\chi_\square)$,
for every prime $q$, $5\le q \le \bound$.
We also recall that \eqref{Dirichlet-class-number} implies that $L(1,\chi_\square)>0$.

 \medskip 
\textbf{Computed data}.  
The lists of values of $L(1,\chi_\square)$, $c_{1}(q)$, $c_{2}(q)$, $c_{3}(q)$, $c_{4}(q)$ 
the corresponding upper bounds for $\beta$, the Littlewood bounds,
the list of cases in which Joshi's inequalities hold, and the values of class number of $\Q(\sqrt{-q})$,
all of them computed for the odd primes $q \le 10^7$,
are available online at the address \indirizzoweb.

\begin{figure}[ht] 
\includegraphics[scale=0.45,angle=0]{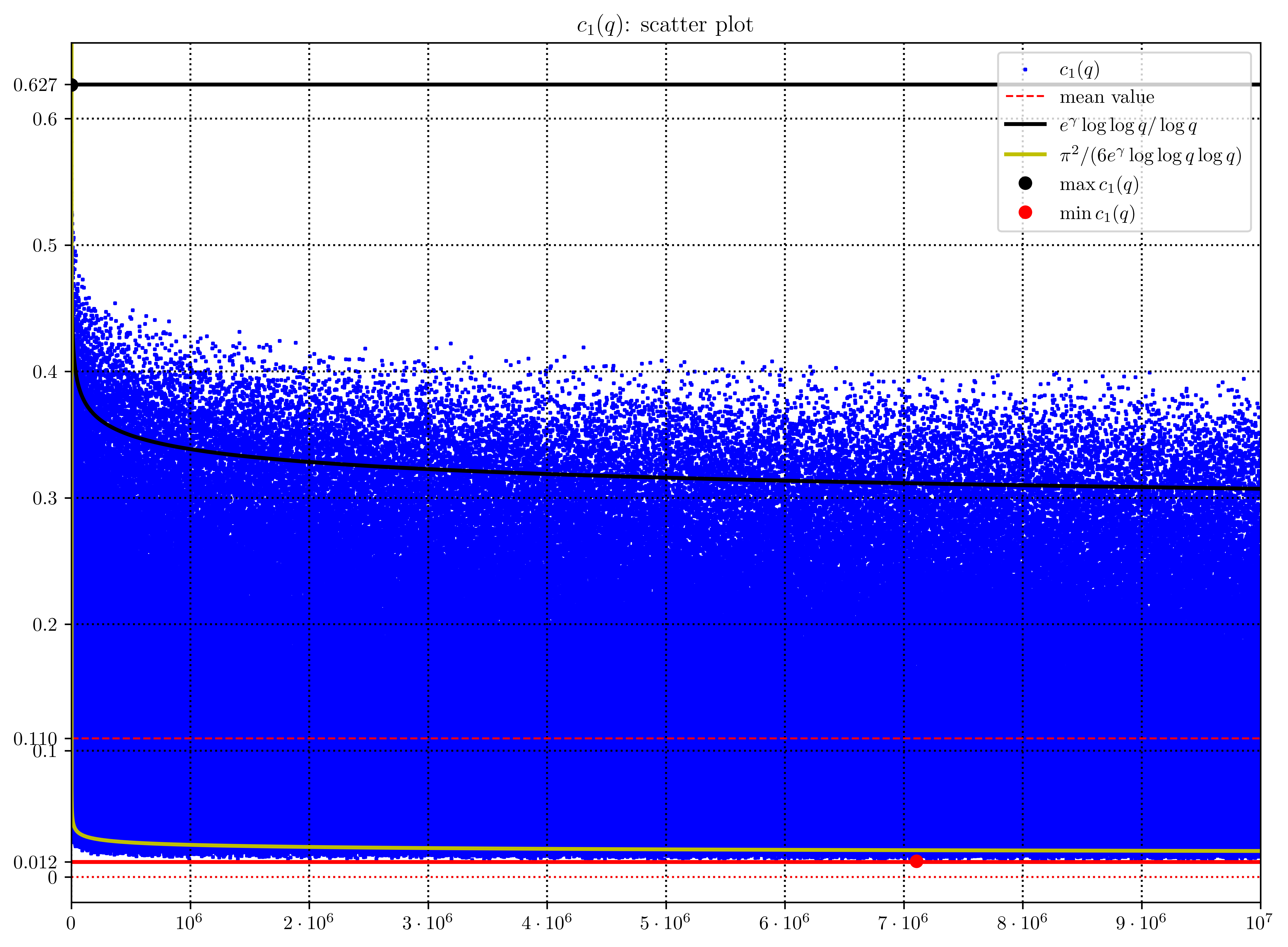}  
\caption{\small{The values of $c_1(q)$, $q$ prime, $3\le q\le \bound$, in Theorem \ref{Thm-beta}.
$ \min c_{1}(q) = \constonemin\dotsc$ attained at 
$q=\constoneminattained$; 
$ \max c_{1}(q) =  \constonemax\dotsc$ attained at 
$q=\constonemaxattained$.
The black straight line corresponds to $0.627$;  the red one  to $0.012$.
The black line corresponds to the first Joshi bound, see \eqref{Joshi-first},
the yellow one corresponds to the second Joshi bound, see \eqref{Joshi-second}.
The red dashed line corresponds to the mean value.
 }}
\label{fig-thm1-c1} 
\end{figure}

\begin{figure}[ht] 
\includegraphics[scale=0.5,angle=0]{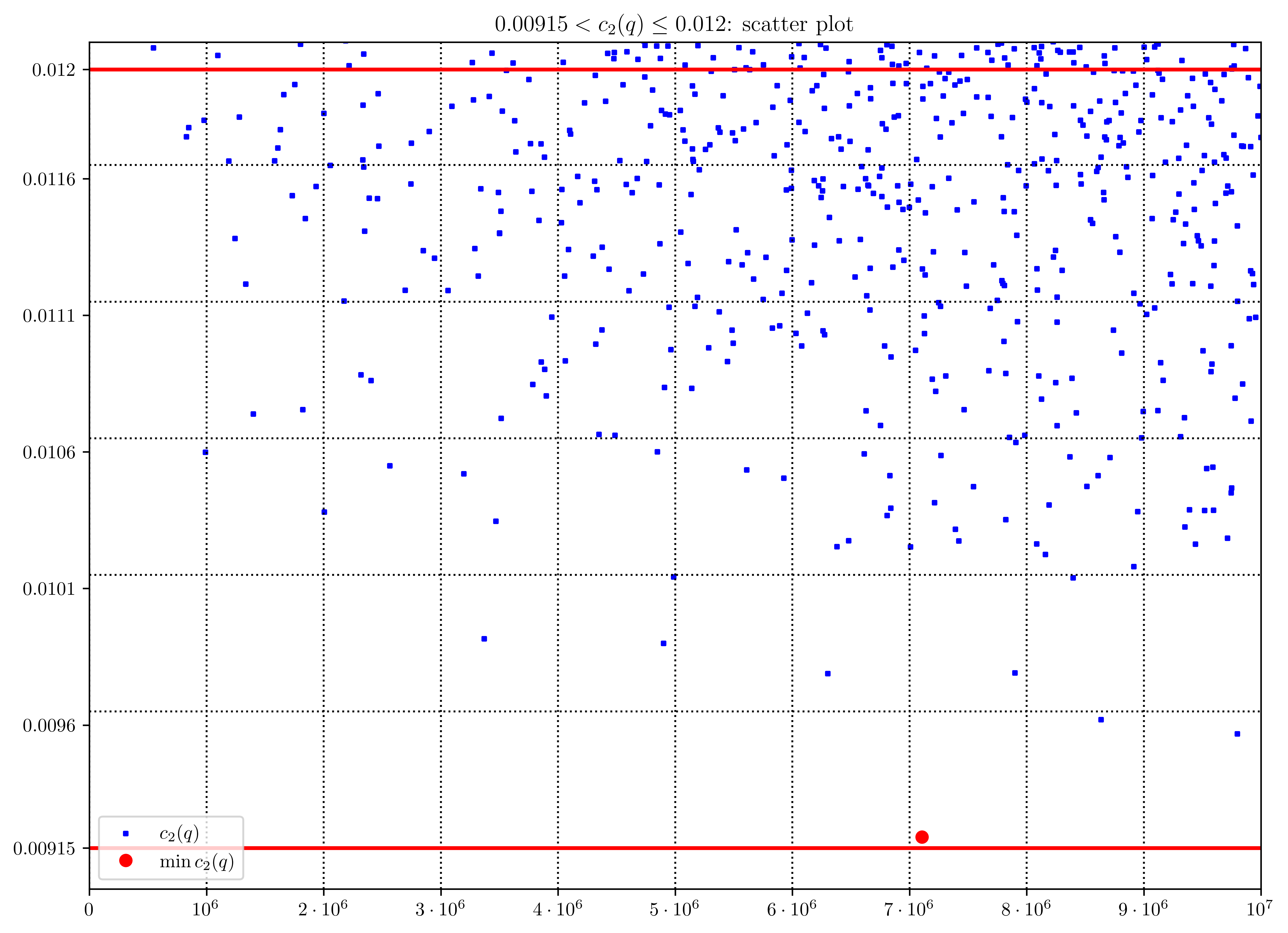}  
\caption{\small{The values of $0.00915<c_2(q) \le 0.012$, $q$ prime, $3\le q\le \bound$, in Theorem \ref{Thm-beta}.
$ \min c_{2}(q) = \consttwomin\dotsc$ attained at 
$q=\consttwominattained$.
The red lines correspond to $0.00915$ and $0.012$.
 }}
\label{fig-thm1-c2-bounded} 
 \end{figure}

\begin{figure}[ht] 
\includegraphics[scale=0.75,angle=0]{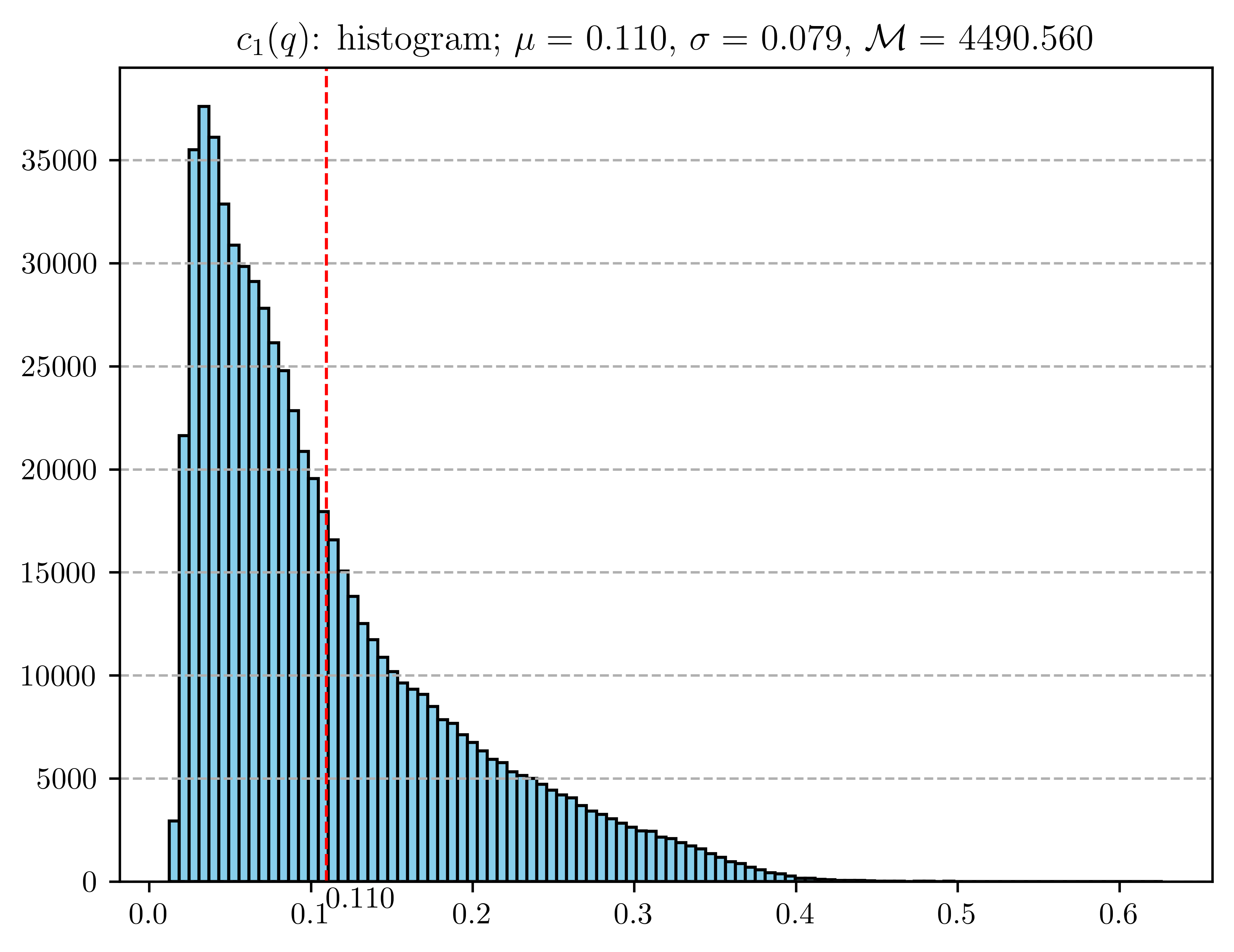}  
\caption{\small{Histogram about the values of $c_1(q)$, $q$ prime, $3\le q\le \bound$, in Theorem \ref{Thm-beta}.
Intervals length $:= I = 0.0067570100\dotsc$;
number of primes $3\le q\le \bound := {\mathcal P} = 664578$;
mass $:=  \mathcal{M} =  I \cdot {\mathcal P}$;
mean $:= \mu = 0.1096877373\dotsc$; 
standard deviation $:= \sigma =0.0788767546\dotsc$
The red dashed line corresponds to the mean value.
 }}
\label{fig-thm1-c1-histo} 
\end{figure}

\begin{figure}[ht] 
\includegraphics[scale=0.75,angle=0]{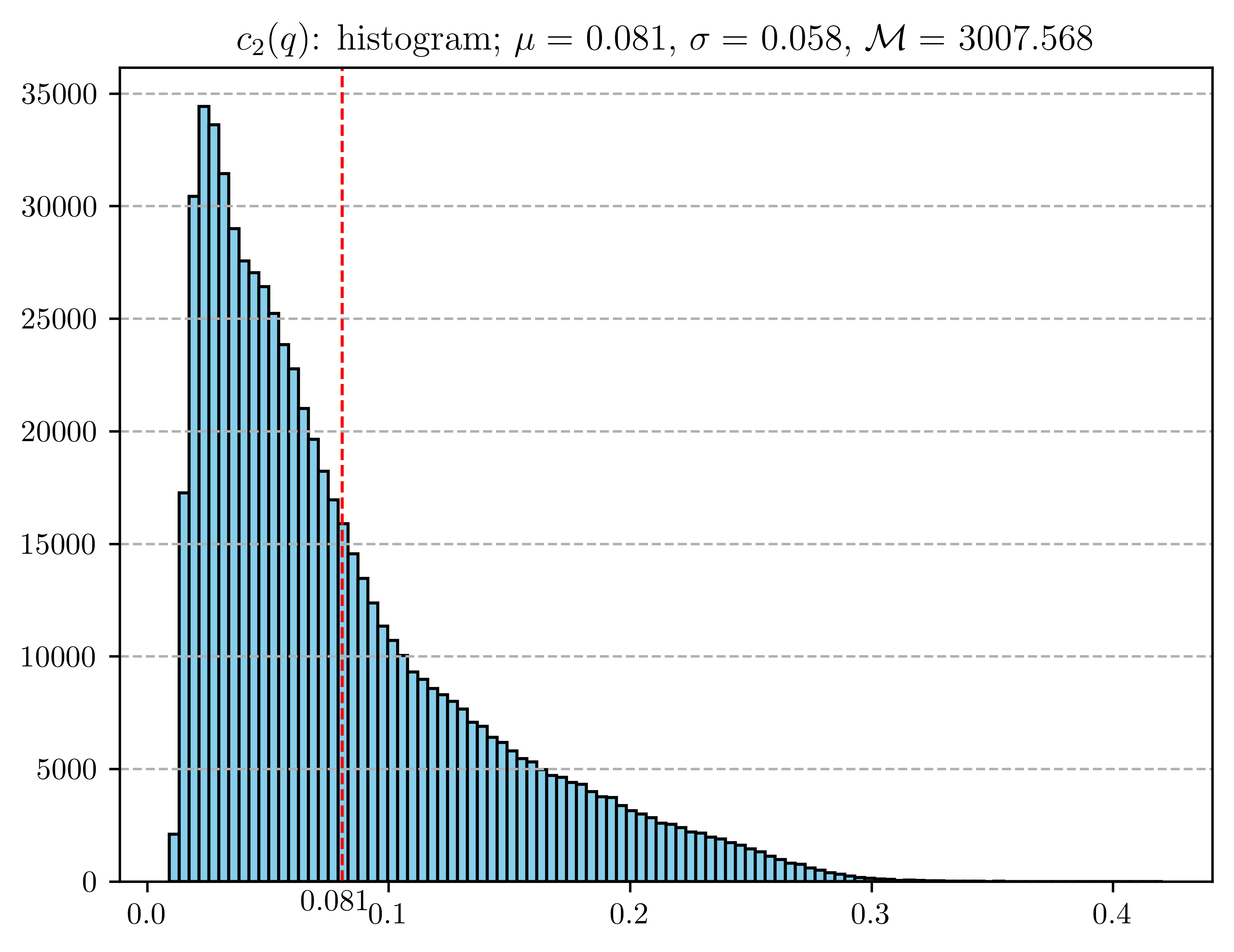}  
\caption{\small{Histogram about the values of $c_2(q)$, $q$ prime, $3\le q\le \bound$, in Theorem \ref{Thm-beta}.
Intervals length $:= I =  0.0045148083\dotsc$;
number of primes $3\le q\le \bound := {\mathcal P} = 664578$;
mass $:=  \mathcal{M} =  I \cdot {\mathcal P}$;
mean $:= \mu = 0.0807233707\dotsc$; 
standard deviation $:= \sigma = 0.0580352818\dotsc$
The red dashed line corresponds to the mean value.
}}
\label{fig-thm1-c2-histo} 
\end{figure}

 \medskip 
\textbf{Tables and Figures}.  
We provide here  the scatter plots of $c_{1}(q)$ and $c_{2}(q)$, 
see Figures \ref{fig-thm1-c1}-\ref{fig-thm1-c2-bounded}.
We also inserted the histograms obtained with the same values,
see Figures \ref{fig-thm1-c1-histo}-\ref{fig-thm1-c2-histo}.
In Figures \ref{fig-ULI}-\ref{fig-LLI} we inserted two scatter plots about $\ULI$ and $\LLI$
and in Figures \ref{fig-ULI-histo}-\ref{fig-LLI-histo} we presented the corresponding
histograms.
The values of $L(1,\chi_\square)$, $c_{1}(q)$, $c_{2}(q)$, $c_{3}(q)$, $c_{4}(q)$ 
and the computed upper bounds for $\beta$
for $3\le q \le \boundtable$, $q$ prime, are  collected in Tables  \ref{TableL}-\ref{Tablec3c4}.

 \medskip 
\textbf{Outline}.  
The paper is organised as follows:
in Section \ref{proof-Thm1} we will prove Theorem \ref{Thm-beta};
a part of its proof is based on the computation described 
in Section \ref{chi-Bernoulli-method}.
In particular, in this section we will see how to compute 
$L(1,\chi_\square)$   using the values 
of Euler's $\Gamma$-function and the Fast Fourier Transform algorithm.
After the bibliography 
we will insert the tables and figures previously mentioned. 
 
 \medskip 
\textbf{Acknowledgements}.  
The harder part of the computation
was performed on the cluster of machines  of the Dipartimento di Matematica
``Tullio Levi-Civita'' of the University of 
Padova, see \clusteraddress.
I wish to thank the IT-manager, Luca Righi,
for his help in organising and handling the actual computation.
I would also like to thank Pieter Moree 
for having read a preliminary version of this paper.
Finally, I would like to thank the reviewer for his/her suggestions as he/she
allowed me to improve the presentation of my work.

\section{Proof of Theorem \ref{Thm-beta}}
\label{proof-Thm1}
The first result we need is the one that follows from the 
computations we will describe in the next section.
\begin{Lemma}
\label{L1-estim-computed}
Let $q$ be an odd prime, $q\le \bound$, let $\chi_\square$ be 
the quadratic Dirichlet character mod $q$.
We have that 
\[
d_1 \log q<  L(1,\chi_\square)   < d_2 \log q,
\]
where $d_1 = 0.0124862668\dotsc$ and  $d_2 = 0.6267599041\dotsc$,
attained at  $q=\constoneminattained$ ($\chi_\square$ is even),  respectively $q = 23$
($\chi_\square$ is odd).
 \end{Lemma}

The complete list of values for $L(1,\chi_\square)$
for each odd prime $q \le 10^7$ is available
online at the web address mentioned in the Introduction.
Their values for $3\le q \le \boundtable$ are collected in Table \ref{TableL}.

We will explain now how  to obtain the estimates on the Landau-Siegel zero $\beta$
from the ones on $L(1,\chi_\square)$, where $L(s,\chi_\square)\in\LL$, as defined in \eqref{L-set-def}.
We will need the following lemmas.
The first one is an explicit version of the P\'olya-Vinogradov estimate
proved by Lapkova \cite[Lemma 3A]{Lapkova2018} in 2018.
\begin{Lemma}
\label{Lapkova}
Let $\chi$ be a non-principal primitive Dirichlet character mod $q\ge 2$.
Define
\begin{equation}
\label{g-funct-def}
g(q) : = 
\begin{cases}
\frac{2}{\pi^2} + \frac{0.9467}{\log q} + \frac{1.668}{\sqrt{q}\log q}  & \text{if}\ \chi \ \text {is even} \\
\frac{1}{2\pi} + \frac{0.8204}{\log q}+ \frac{1.0286}{\sqrt{q}\log q}&\text{if}\ \chi \ \text {is odd}.
\end{cases}
\end{equation}
Then
\begin{equation}
\label{S-chi-def}
S(\chi,q) : =
\max_{N} \Bigl \vert \sum_{a=q+1}^{N} \chi(a) \Bigr \vert
\le
g(q) \sqrt{q}  \log q .
\end{equation}
\end{Lemma}
Before proceeding further,
we need to define the following quantity:
\begin{equation}
\label{S-def}
S(q) : = \sum_{n=2}^{q} \frac{\log n}{n}.
\end{equation}
The second lemma we need is about a suitable estimate on $L^\prime(\sigma, \chi)$
for $\sigma$ ``close'' to $1$.
\begin{Lemma}
\label{Lprime-estimate}
Let $\sigma$ be such that $1-1/\log q \le \sigma  \le 1$.
Let $\chi$ be a non-principal Dirichlet character mod $q\ge 2$.
Then
\[
\vert L^\prime(\sigma, \chi) \vert 
\le
e
\Bigl(
S(q) + \frac{\log q}{q} S(\chi,q)
\Bigr),
\]
where $S(q)$ is defined in \eqref{S-def} and
 $S(\chi,q)$ is defined in \eqref{S-chi-def}.
\end{Lemma}
\begin{Proof}[Proof of Lemma \ref{Lprime-estimate}]
We argue as in Davenport \cite[p.~96, proof of (11)]{Davenport2000}.
Recalling $L^\prime(\sigma, \chi) = -\sum_{n=2}^{\infty}  \chi(n) n^{-\sigma} \log n$,
we split the sum according to $n\le q$ and $n\ge q+1$.
For the first, using
$n^{-\sigma} \le e/n$ for $n \le q$ and  $\sigma \ge 1-1/\log q$,
we obtain
\[
\Bigl\vert 
\sum_{n=2}^{q} \frac{\chi(n) \log n}{n^\sigma}
\Bigr\vert
\le
e S(q) ,
\]
where $S(q)$ is defined in \eqref{S-def}.
For the second, using the partial summation formula
and $q^{-\sigma} \le e/q$   for  $\sigma \ge 1-1/\log q$,
we have
\[
\Bigl\vert 
\sum_{n=q+1}^{\infty} \frac{\chi(n) \log n}{n^\sigma}
\Bigr\vert
\le
\frac{\log q}{q^\sigma} \max_N \Bigl \vert \sum_{n=q+1}^{N} \chi(n)  \Bigr \vert
\le
e \frac{\log q}{q} S(\chi,q),
\]
where
 $S(\chi,q)$ is defined in \eqref{S-chi-def}.
Lemma \ref{Lprime-estimate} immediately follows.
\end{Proof}

Reasoning as in Davenport, see \cite[p.~95--96]{Davenport2000},
we can now prove Theorem \ref{Thm-beta}.
Let  $1-1/\log q \le \sigma \le 1$ and assume that $L(\beta,\chi_\square)=0$.
By the mean value theorem we obtain
\[
L(1,\chi_\square) - L(\beta,\chi_\square)
= 
(1-\beta ) L^\prime(\sigma, \chi_\square).
\]
Assuming there exists $c_{1}(q)>0$ such that 
\begin{equation}
\label{L1-lower-bound}
 L(1,\chi_\square) > c_{1}(q) \log q
\end{equation}
and using Lemma \ref{Lprime-estimate}, we obtain 
\begin{align*}
c_{1}(q) \log q 
&< 
(1-\beta) \vert L^\prime(\sigma, \chi_\square)\vert
\le
e (1-\beta)  \Bigl(S(q)  + \frac{\log q}{q} S(\chi_\square,q)\Bigr)
\\&
\le
e (1-\beta)  \Bigl(S(q) + \frac{(\log q)^2}{\sqrt{q}} g(q)\Bigr),
\end{align*}
where $g(q)$ is defined in \eqref{g-funct-def}
and in the last step we used Lemma \ref{Lapkova}.
Letting
\begin{equation}
\label{c3-c4-def}
c_{3}(q) := e \frac{S(q) }{(\log q)^2}
\quad
\text{and}
\quad
c_{4}(q) := e \frac{g(q)}{\sqrt{q}},
\end{equation}
where $S(q)$ is defined in \eqref{S-def},
we obtain
\[
c_{1}(q) \log q  <  (1-\beta) (c_{3}(q) + c_{4}(q))(\log q)^2.
\]
Hence
\begin{equation}
\label{c1-c2-function}
\beta < 1 - \frac{c_{2}(q)}{\log q},
\quad
\text{where}
\quad
c_{2}(q) :=\frac{c_{1}(q)}{c_{3}(q) + c_{4}(q)},
\end{equation}
and $c_{1}(q),c_{3}(q), c_{4}(q)$ are respectively defined in 
\eqref{L1-lower-bound}-\eqref{c3-c4-def}.

Using Lemma \ref{L1-estim-computed} we have,
for every odd prime $q\le 10^7$, that
$c_{1}(q)  = d_1$.
Moreover, both $c_{3}(q)$ and $c_{4}(q)$ can be easily computed,
for example using Pari/GP \cite{PARI2022}.\footnote{The complete list of values for $c_{3}$ and $c_{4}$ 
for every odd prime $q\le10^7$ are available
online at the web address mentioned in the Introduction.
Their values for $3\le q \le \boundtable$, $q$ prime, are collected in Table \ref{Tablec3c4}.}
Theorem 1 hence follows. 

\section{Computation of $\vert L(1,\chi)\vert$ and  proof of Theorem \ref{Thm-beta}}
\label{chi-Bernoulli-method} 
Recall that $q$ is an odd prime and let $\chi$ be a primitive  non-principal 
Dirichlet character mod $q$.
Distinguishing between the parity of the Dirichlet characters,
 if $\chi$ is  even  we have, see, \emph{e.g.},
 Cohen \cite[proof of Proposition 10.3.5]{Cohen2007}, that
 \[
 L(1,\chi)
= 
2 \frac{\tau(\chi)}{q}
\sum_{a=1}^{q-1} \overline{\chi}(a)\log\Bigl(\Gamma\bigl(\frac{a}{q}\bigr)\Bigr)
\] 
where the \emph{Gau\ss\ sum} 
$\tau(\chi):= \sum_{a=1}^q \chi(a)\, e(a/q)$, $e(x):=\exp(2\pi i x)$, satisfies 
$\vert \tau(\chi) \vert = \sqrt{q}$. Hence
\begin{equation}
\label{even} 
\vert 
L(1,\chi)
\vert 
= 
\frac{2}{\sqrt{q}}\
\Bigl\vert 
\sum_{a=1}^{q-1} \overline{\chi}(a)
\log\Bigl(\Gamma\bigl(\frac{a}{q}\bigr)\Bigr)\Bigr\vert
\quad 
(\chi\ \text{even}).
\end{equation}
Moreover, if $\chi$ is an odd character,  we have, see, \emph{e.g.},
Cohen \cite[Corollary 10.3.2]{Cohen2007}, that
\[
L(1,\chi)
=
-   \frac{\pi\tau(\chi)} {q} \sum_{a=1}^{q-1}  \frac{a\overline{\chi}(a)}{q} 
\] 
and hence we obtain
\begin{equation}
\label{odd}
\vert 
L(1,\chi)
\vert 
=
\frac{\pi}{\sqrt{q}} \ \Bigl\vert \sum_{a=1}^{q-1}  \frac{a\overline{\chi}(a)}{q}  \Bigr\vert
\quad 
(\chi\ \text{odd}).
\end{equation}   

In both the equations \eqref{even}-\eqref{odd} we can embed a \emph{decimation 
in frequency strategy} in the Fast Fourier Transform (FFT) algorithm used to perform the sum over $a$,
see Section \ref{fft-approach}.
The needed set of Gamma-function  values can be computed with a precision of $n$ binary digits with a cost of 
$\Odi{qn}$ floating point products,  see \cite[Section 3]{Languasco2021}, plus the cost of computing $(q-1)/2$ 
values of the logarithm function.
Hence, recalling also that the computational cost of the FFT algorithm 
of length $q$ is $\Odi{q\log q}$ floating point products, 
the total cost for computing $\vert L(1,\chi)\vert$ with a precision of $n$ binary digits
is then  $\Odi{q(n+\log q)}$ floating point products.
Currently, this is
the fastest algorithm  to compute $\vert L(1,\chi)\vert$. 

We now proceed to describe our computational strategy. 

 \subsection{Computations using Pari/GP (slower, more decimal digits available).}
In practice we first computed a few values of $ L(1,\chi_\square)$ using Pari/GP, v.~2.15.1, 
since it has  the ability to  generate the Dirichlet $L$-functions (and  many other $L$-functions).
This can be done with few instructions of the gp scripting language.
Such a computation has a linear cost in the number of calls
of the {\tt lfun} function of Pari/GP and it is, at least
on our Dell Optiplex desktop machine, slower than 
 using \eqref{even}-\eqref{odd}.
We used the  {\tt lfun}-approach to get the values of $L(1,\chi_\square)$  for every $q$ prime, $3\le q\le \boundtable$, 
with a precision of  $38$ decimal digits (see Table \ref{TableL})
in less than $1$ second of computation time. 
The machine we used was a Dell OptiPlex-3050, equipped with an 
Intel i5-7500 processor, 3.40GHz,  32 GB of RAM and running Ubuntu 22.04.1 LTS.

In fact, since just one \texttt{lfun}-call is needed for each $q$, we extended
such a computation to every $q\leq 10^6$ with a precision of 38 decimal digits. 
That required about $3$ hours and $12$ minutes of computing time on the machine previously mentioned.
   
\subsection{Building the FFT approach}
\label{fft-approach} 
 As $q$ becomes large, the time spent in summing over $a$
 in both \eqref{even} and \eqref{odd} 
dominates the overall computational cost.  So we implemented
the use of the FFT 
by using  the {\tt fftw} \cite{FFTW} library in our   C programs.   
We will explain now how to do so.

 In both   \eqref{even} and \eqref{odd} we remark that,
since $q$ is prime, it is enough to get $g$, a primitive root of $q$,
and $\chi_1$, the Dirichlet character mod $q$ given by 
 $\chi_1(g) = e^{2\pi i/(q-1)}$, to see that the set of the non-principal characters
 mod $q$ is $\{\chi_1^j \colon j=1,\dotsc,q-2\}$.
It is well known that the problem of finding a primitive root $g$ of $q$ is a 
computationally hard one, see, e.g., Shoup \cite[Chapters 11.1-11.4]{Shoup2005},
but, for each fixed prime $q$, we need to find it just once. In the applications,
for each involved prime $q$ we can save such a $g$ and reuse
it every time we have to work again mod $q$.

Hence, if, for every $k\in \{0,\dotsc,q-2\}$, we denote $g^k\equiv a_k \bmod{q}$, 
every summation in  \eqref{even}  and \eqref{odd}
is  of the type $\sum_{k=0}^{q-2}  e^{-2\pi i j k /(q-1)} f(a_k/q)$, 
where $j\in\{1,\dotsc,q-2\}$ and $f$ is either $\log\Gamma$ or the identity function
according to the parity of the quadratic Dirichlet character $\chi_\square$.
 
The quadratic Dirichlet character $\chi_\square$ has order $2$;
so, in the previous notation, it means that $j = (q-1)/2$ or, in other words,  $\chi_\square = \chi_1^{(q-1)/2}$.
Hence  $\chi_\square (g) = e^{\pi i} = -1$ and $\chi_\square (-1) = \chi_\square (g^{(q-1)/2})
= (-1)^{(q-1)/2}$; so $\chi_\square$ is an even character if and only if $q \equiv 1 \bmod 4$.
This fact, not surprisingly,  means that
the parity of $\chi_\square$ depends on the value of the Legendre symbol $(-1\vert q)$.
  
Using the setting previously described, we can now write that the quantities in \eqref{even}-\eqref{odd}  are
the Discrete Fourier Transform (DFT)  of the sequence $\{ f(a_k/q)\colon k=0,\dotsc,q-2\}$.
 
Even if in this application we are interested in just one value (the one corresponding
to the quadratic character), the use of the FFT algorithm  is preferable,
due to its speediness, at least whenever the main parameter (in this case, $q$)
becomes large. At the end of the FFT computation, it will be enough to
pick up the value corresponding to $j = (q-1)/2$ in the output sequence.
In fact, we obtained such values as a part of a wider computation
involving several other interesting quantities, for example the ones
related to the Euler-Kronecker constants of the cyclotomic field 
and of some of its subfields (for some application, see, e.g., \cite{CiolanLM2023},
\cite{Languasco2021a} and \cite{LanguascoR2021}).

To further speed-up the computation and to reduce its memory usage, we used 
the \emph{decimation in frequency} strategy. Let $f$ be a
function, $f\colon (0,1) \to \R$. Arguing as in \cite[Section 2.2]{Languasco2021}, 
letting $e(x):=e^{2\pi i x}$, $\overline{q}=(q-1)/2$,    
for every $j=0,\dotsc, q-2$, $j=2t+\ell$, $\ell\in\{0,1\}$ and $t\in \Z$, we have that  
\[
 \sum_{k=0}^{q-2}   e\Bigl(\frac{- j k}{q-1}\Bigr)  f \Bigl(\frac{a_k}{q}\Bigr)
 =
 \sum_{k=0}^{\overline{q}-1}  
 e\Bigl(\frac{- t k}{\overline{q}}\Bigr)  
    e\Bigl(\frac{- \ell k}{q-1}\Bigr)    
 \Bigl(
 f\Bigl(\frac{a_k}{q}\Bigr) 
 +
 (-1)^{\ell} 
 f \Bigl(\frac{a_{k+\overline{q}}}{q}\Bigr)
 \Bigr),
\]
 where $t=0,\dotsc, \overline{q}-1$. 
 Letting
 \begin {equation} 
\label{bk-ck-def} 
b_k :=
f\Bigl(\frac{a_k}{q}\Bigr) +  f \Bigl(\frac{a_{k+\overline{q}}}{q}\Bigr)   
\quad
\textrm{and}
\quad
c_k :=  
e\Bigl(-\frac{k}{q-1}\Bigr)   
\Bigl(  f\Bigl(\frac{a_k}{q}\Bigr) -  f \Bigl(\frac{a_{k+\overline{q}}}{q}\Bigr)  \Bigr),
\end{equation}
we can rewrite the previous formula  (recall that 
$j=2t+\ell$, $\ell\in\{0,1\}$ and $t=0,\dotsc, \overline{q}-1$) as
\begin {equation} 
\label{DIF} 
\sum_{k=0}^{q-2}   e\Bigl(\frac{- j k}{q-1}\Bigr)  f \Bigl(\frac{a_k}{q}\Bigr)
 =
\begin{cases}
\sum\limits_{k=0}^{\overline{q}-1}     e\bigl(-\frac{t k}{\overline{q}}\bigr) b_k 
& \textrm{if} \ \ell =0\\
\sum\limits_{k=0}^{\overline{q}-1}     e\bigl(-\frac{t k}{\overline{q}}\bigr)  c_k  
& \textrm{if} \ \ell =1.\\
\end{cases}
\end{equation}
Since we just need the sum  over the odd Dirichlet characters for $f(x)=x$ 
and over the even Dirichlet characters for $f(x)=\log \Gamma(x)$,
in this way we can evaluate  an FFT of length $\overline{q}=(q-1)/2$, instead of $q-1$,
applied on a suitably modified sequence according to \eqref{bk-ck-def}-\eqref{DIF}.
Clearly this represents a gain in both the speed and the memory usage
in running the actual computer program. 

In the case $f(x)= \log \Gamma(x)$ we can simplify the form of 
$b_k = \log \Gamma(a_k/q) + \log \Gamma(a_{k+\overline{q}}/q) $, where 
$\overline{q}=(q-1)/2$ and $k=0,\dotsc, \overline{q}-1$,
in the following way. Recalling that $\langle g \rangle = \Z^*_q$, 
$a_k \equiv g^k \bmod q$ and   $g^{\overline{q}} \equiv q-1 \bmod{q}$,
we can write that  $ a_{k+\overline{q}}  \equiv  q-a_{k} \bmod{q}$
and hence
\(
\log \Gamma(a_{k+\overline{q}}/q)  = \log \Gamma ((q-a_{k})/q)
=
\log \Gamma (1- a_{k}/q).
\)
Using the well-known \emph{reflection formula} 
$\Gamma(x) \Gamma(1-x)  = \pi / \sin(\pi x)$,  
we obtain 
\begin{align*}
%
\log \Gamma\Bigl(\frac{a_{k}}{q}\Bigr)
+
\log \Gamma\Bigl(\frac{a_{k+\overline{q}}}{q}\Bigr)
=
\log \Gamma\Bigl(\frac{a_{k}}{q}\Bigr)
+
\log \Gamma\Bigl(1-\frac{a_{k}}{q}\Bigr)
&=
 \log \pi - \log \Bigl(\sin\bigl( \frac{\pi a_k}{q}\bigr)\Bigr) ,
\end{align*}
for every $k=0,\dotsc, \overline{q}-1$. Inserting the last relation in the definition of $b_k$  in \eqref{bk-ck-def} 
and remarking that, by the orthogonality of the Dirichlet characters, the final contribution of the constant term $\log \pi$ will be zero,
we can replace in the actual computation  the  Gamma function  with the $\log(\sin (\cdot))$ one.
Since in our application we will have $a/q\in(0,1)$,  
we used our own alternative implementation of $\log \Gamma(x)$, $x\in (0,1)$, 
see \cite[Section 4]{Languasco2023a}, because in this way we have a further 
gain in the execution speed.\footnote{This is due to the fact that we can use a series in which half of the coefficients
are equal to zero and we can also accurately control how many summands are
needed to computing such a sequence of $\log\Gamma$-values up to the required accuracy.
See the reflection formulae for  $\log \Gamma(x)$
in Proposition 3 of \cite{Languasco2023a}.}

In the case $f(x)=x$, it is easier to obtain a simplified  form of $c_k$  as defined in \eqref{bk-ck-def}. 
Using again $\langle g \rangle = \Z^*_q$, 
$a_k \equiv g^k \bmod q$ and   $g^{\overline{q}} \equiv q-1 \bmod{q}$,
 we can write that $ a_{k+\overline{q}}  \equiv  q-a_{k} \bmod{q}$;  hence
\(
a_k  -    a_{k+\overline{q}}  
= 
a_k -(q-a_{k})  
=
2a_k  -q.
\)
Summarising,  for every $k=0,\dotsc, \overline{q}-1$, $\overline{q}=(q-1)/2$, we obtain 
\[
c_k=   e\Bigl(-\frac{k}{q-1}\Bigr)\Bigl(2\frac{a_k}{q} -1\Bigr).
\]
  
\subsection{Computations  summing over $a$ via FFT (much faster, less decimal digits available).}
Using the setting explained in the previous subsection,  we were able to compute, 
using the \emph{long double precision} (80 bits)
of the C programming language, the values of $ L(1,\chi_\square) $ 
for every odd prime $q \le  \bound$ and we provide here 
the scatter plots of such values normalised using $\log q$ (in other words,
the values of $c_1(q)=  L(1,\chi_\square)  / \log q$, see \eqref{L1-lower-bound}), 
and the corresponding values for $c_2(q)$
as defined in \eqref{c1-c2-function},
see Figures \ref{fig-thm1-c1}-\ref{fig-thm1-c2-bounded}.
We also inserted the histograms obtained with the same values,
see Figures \ref{fig-thm1-c1-histo} and \ref{fig-thm1-c2-histo}.
The data were obtained  in about 
$4$ days  of global computation time on 
the cluster of machines  of the  Dipartimento di Matematica
``Tullio Levi-Civita'' of the University of 
Padova, see \clusteraddress.
The actual FFTs were performed using the FFTW \cite{FFTW} software library.

\subsection{FFT accuracy estimate: $q=9999991$}
According to Schatzman \cite[\S~3.4, p.~1159-1160]{Schatzman1996},  
the root mean square relative error in the FFT is bounded by  
\begin{equation}
\label{Delta-FFT-9999991}
\Delta = \Delta(N, \eps) :=0.6 \eps \sqrt{\log_2 N} ,
\end{equation} 
where $\eps$ is the machine epsilon,
$N$ is the length of the transform and $\log_2 x $ denotes
the logarithm in base $2$ of $x$. 
According to the IEEE 754-2008 specification, we can 
set $\eps=2^{-64}$ for the \emph{long double precision} 
of the C programming language. 
So for the largest prime less than $10^7$, $q=9999991$, 
letting $N=(q-1)/2$, we get that $\Delta<1.54 \cdot 10^{-19}$.
To evaluate the euclidean norm of the error we have then to multiply $\Delta$ and the
euclidean norms of the  sequences
\begin{alignat*}{3}
x_k&:= 2\frac{a_k}{q}-1, \quad
&y_k&:= \log\Gamma\Bigl(\frac{a_k}{q}\Bigr) + \log\Gamma\Bigl(1-\frac{a_k}{q}\Bigr) - \log \pi, 
\end{alignat*}
where  $a_k  = g^k \bmod q$, $\langle q \rangle =\Z_q^*$.
A straightforward computation gives 
\[
\Vert x_k \Vert_2 = \Bigl(   \frac{(q-1)(q-2)}{6q} \Bigr)^{1/2}
=1290.99367\dotsc.
\]
Hence, recalling that $\Vert \cdot \Vert_{\infty} \le \Vert \cdot \Vert_{2}$,
for this sequence we can estimate that the maximal error in its FFT-computation
is bounded by $1.99 \cdot 10^{-16}$ (long double precision case). 
Unfortunately, no closed formulas for the euclidean norm of  $y_k$ are known 
but, using  $\Vert \cdot \Vert_{\infty} \le \Vert \cdot \Vert_{2} \le \sqrt{N} \Vert \cdot \Vert_{\infty}$
and 
\begin{align*}
\Vert y_k \Vert_{\infty} &=  - \log \sin(\frac{\pi}{q}) =  14.97336\dotsc,
\end{align*}
that can be obtained using straightforward computations,
we have that  the error in its  FFT-computation is
$< 5.14 \cdot 10^{-15}$.

We also estimated \emph{in practice} the accuracy in 
the actual computations  using the 
FFTW software library by evaluating at   
run-time the quantities
\(
\mathcal{E}_{j}(u_k) : = 
\Vert \mathcal{F}^{-1}(\mathcal{F}(u_k)) - u_k \Vert_{j},
\)
$j\in\{2,\infty\}$,
where $u_k$ runs over the sequences $\{x_k,y_k\}$, 
$\mathcal{F}(\cdot)$ is the Fast Fourier Transform 
and $\mathcal{F}^{-1}(\cdot)$ is its inverse transform.
Theoretically we have that  $\mathcal{E}_{j}(u_k) =0$;
moreover, assuming that the root mean square relative error in the FFT is bounded by  
$\Delta>0$, it is easy to obtain
\begin{equation}
\label{back-forth-estim-9999991}
\mathcal{E}_{2}(u_k) < \Delta(2+\Delta) \Vert u_k \Vert_{2}
\quad
\textrm{and}
\quad
\mathcal{E}_{\infty}(u_k) < \Delta(2+\Delta)  \sqrt{N}  \Vert u_k \Vert_{\infty}.
\end{equation}
For  $q=9999991$, $N=(q-1)/2$ and $\eps = 2^{-64}$ in \eqref{Delta-FFT-9999991},
we get that $\Delta(2+\Delta)< 3.07 \cdot 10^{-19}$
and, using again the previous norm-values, 
we also obtain that $\mathcal{E}_{\infty}(u_k)< 1.03 \cdot 10^{-14}$,
where $u_k$  runs over the sequences $\{x_k,y_k\}$.
Moreover, the actual computations using FFTW gave the following results:  
\begin{align*}
\frac{\mathcal{E}_2(x_k)}{\Vert x_k \Vert_2} < 3.57\cdot 10^{-19}, \quad
\frac{\mathcal{E}_2(y_k)}{\Vert y_k \Vert_2} < 3.06\cdot 10^{-19}, \quad
\end{align*}
in agreement with the first part of \eqref{back-forth-estim-9999991},
and
\begin{align*}
\mathcal{E}_\infty(x_k) <  8.25 \cdot 10^{-19}, \quad
\mathcal{E}_\infty(y_k) <  3.53 \cdot 10^{-18}, \quad
\end{align*}
which are in fact much better than  the second part of \eqref{back-forth-estim-9999991}.
Summarising, we can conclude that at least ten decimal 
digits of our final results are correct.
If necessary,  we can repeat the FFT-step 
using the 
\emph{quadruple precision} ($128$ bits) numerical type:
the final results will be more accurate since   in this case we can use
$\eps = 2^{-113}$ in \eqref{Delta-FFT-9999991} but 
the actual computing time will be much longer.

\subsection{Numerical values in the statement of Theorem \ref{Thm-beta}}

The results obtained as described before in this section were then collected in some comma-separated values (csv) files
together with the ones for $S(q)$ and $S(q,\chi)$ obtained with some 
straighforward Pari/GP scripts. 

A suitable program written in  \texttt{python},~v.3.11.3,
and using the package  \texttt{pandas},~v.2.0.0,
performed the analysis to obtain the  values mentioned in the statement
of Theorem \ref{Thm-beta} and in the Introduction of this paper. The same program also 
gave us the plots we will show in the next section;
they were obtained using  the package \texttt{matplotlib},~v.3.7.1.

The Pari/GP scripts, the C and python programs used and the computational 
results\footnote{In such files the least significant digit of each result might be rounded by the python or the Pari/GP printing/saving routines.} obtained are available at the web address mentioned in the Introduction.

\vskip 0.25cm
\noindent
Alessandro Languasco
Universit\`a di Padova,
Dipartimento di Matematica,
``Tullio Levi-Civita'',
Via Trieste 63,
35121 Padova, Italy.   
{\it e-mail}: alessandro.languasco@unipd.it 

\begin{table}[htp]
\scalebox{0.75}{
\begin{tabular}{|r|c|}
\hline
$q$ & $ L(1,\chi_\square)  $ \\ \hline
3& 0.60459978807807261686\\ 
5& 0.43040894096400403888\\ 
7& 1.18741041172372594878\\ 
11& 0.94722582509948293642\\ 
13& 0.66273539107184558971\\ 
17& 1.01608483384284075219\\ 
19& 0.72073078414566794539\\ 
23& 1.96520205410785916590\\ 
29& 0.61176628956230686982\\ 
31& 1.69274009217927609028\\ 
37& 0.81929216872543187792\\ 
41& 1.29909306157509892164\\ 
43& 0.47908838823985721176\\ 
47& 2.29124192852861593669\\ 
53& 0.54002494510255821476\\ 
59& 1.22700157894863547525\\ 
61& 0.93831019824883536614\\ 
67& 0.38380662888291551638\\ 
71& 2.60986917715784586434\\ 
73& 1.79463648373514640544\\ 
79& 1.76728394209597891666\\ 
83& 1.03450377843099381393\\ 
89& 1.46444140226401940465\\ 
97& 1.89349532374748011332\\ 
101& 0.59666866801751914350\\ 
103& 1.54775161082393862005\\ 
107& 0.91112767558291391313\\ 
109& 1.06597159422690850687\\ 
113& 1.38235170906160842756\\ 
127& 1.39385634554837260497\\ 
131& 1.37241112290895862735\\ 
137& 1.39381341332016437677\\ 
139& 0.79939923310163592599\\ 
149& 0.67359583330192465588\\ 
151& 1.78961429055614439182\\ 
157& 0.85575891033757795011\\ 
163& 0.24606852755296024389\\ 
167& 2.67414112075698568624\\ 
173& 0.39091084118993365939\\ 
179& 1.17406829821234371721\\ 
181& 1.06647233219024291485\\ 
191& 2.95512966360404799352\\ 
193& 2.17043401934037005992\\ 
197& 0.47500088853260044973\\ 
199& 2.00431438732986139006\\ 
211& 0.64882847253825363656\\ 
223& 1.47263623115067414101\\ 
227& 1.04257413986938056150\\ 
229& 1.07546851605294369728\\ 
233& 1.40761515138248090751\\ 
239& 3.04819103378239805449\\ 
241& 2.41835638390008713223\\ 
251& 1.38806898983133714208\\ 
257& 1.29748495883097760616\\ 
263& 2.51834572398050237072\\ 
269& 0.62189322157495099927\\ 
\hline
\end{tabular}
\begin{tabular}{|r|c|}
\hline
$q$ & $ L(1,\chi_\square)  $ \\ \hline
271& 2.09921979227515633845\\ 
277& 0.94551517859218221604\\ 
281& 1.73837686726620650326\\ 
283& 0.56024489726455251491\\ 
293& 0.33143839429193379380\\ 
307& 0.53790048971882632435\\ 
311& 3.38472414241331604469\\ 
313& 2.18765235235254611623\\ 
317& 0.50422804930237370784\\ 
331& 0.51803264722850692708\\ 
337& 2.33496354618023864501\\ 
347& 0.84324765016103968640\\ 
349& 1.05143274302898396424\\ 
353& 1.26326462707244921597\\ 
359& 3.15033145392974552469\\ 
367& 1.47590821481251243553\\ 
373& 0.95620222669031687010\\ 
379& 0.48411832547443343593\\ 
383& 2.72897405712794362791\\ 
389& 0.79595304655678003747\\ 
397& 0.81759772693728272404\\ 
401& 1.84245030922576326425\\ 
409& 2.58450736883986784592\\ 
419& 1.38129159965066196073\\ 
421& 1.26771733328139755787\\ 
431& 3.17782906126902242265\\ 
433& 2.24855924140863418396\\ 
439& 2.24910054927309989676\\ 
443& 0.74630785716666530102\\ 
449& 1.86439354511048221777\\ 
457& 2.38525593186943368168\\ 
461& 0.54957192848180195934\\ 
463& 1.02201534654906780970\\ 
467& 1.01762899434701362464\\ 
479& 3.58857580472017716180\\ 
487& 0.99651406339751772517\\ 
491& 1.27600282553383958283\\ 
499& 0.42191100603520136968\\ 
503& 2.94161055306196188762\\ 
509& 0.60545067376348094520\\ 
521& 1.69667214053319208099\\ 
523& 0.68686127642631349436\\ 
541& 1.21666267245791056249\\ 
547& 0.40297440642054359518\\ 
557& 0.46302148873976553639\\ 
563& 1.19162110015637966561\\ 
569& 1.88476684215394745733\\ 
571& 0.65735780364210046723\\ 
577& 2.25649569369048617219\\ 
587& 0.90767184035159689867\\ 
593& 1.14973060344194813822\\ 
599& 3.20904989883664076353\\ 
601& 2.71357102171217499608\\ 
607& 1.65767305982970466703\\ 
613& 0.92899988878750402864\\ 
617& 1.46721777647419322522\\ 
\hline
\end{tabular}
\begin{tabular}{|r|c|}
\hline
$q$ & $ L(1,\chi_\square)  $ \\ \hline
619& 0.63135634954328056445\\ 
631& 1.62584277488920558388\\ 
641& 1.97514059108196598356\\ 
643& 0.37167696054672361507\\ 
647& 2.84070128925957940909\\ 
653& 0.58035647957020970787\\ 
659& 1.34616982230635910037\\ 
661& 1.11999325199804251383\\ 
673& 2.48337613840124770910\\ 
677& 0.30374567568477932175\\ 
683& 0.60104851043764983729\\ 
691& 0.59755908466288596373\\ 
701& 0.76048659630628182151\\ 
709& 1.30788527015954699221\\ 
719& 3.63201071708524524489\\ 
727& 1.51469788698487771379\\ 
733& 0.73071044851973008699\\ 
739& 0.57782676937727361509\\ 
743& 2.42033098436211130745\\ 
751& 1.71957516129605273827\\ 
757& 1.07750048093693843269\\ 
761& 1.60466139385922576087\\ 
769& 2.74259278564652916003\\ 
773& 0.35496510586623235281\\ 
787& 0.55992841939873066257\\ 
797& 0.41835782349411989520\\ 
809& 1.93293173914226876225\\ 
811& 0.77221367037683702464\\ 
821& 0.67651235685636589623\\ 
823& 0.98558132212062919454\\ 
827& 0.76470716246213692316\\ 
829& 1.19815159165959413711\\ 
839& 3.57917416660357193379\\ 
853& 0.69994286174949083151\\ 
857& 1.13428247579320960158\\ 
859& 0.75032830660913181702\\ 
863& 2.24576220996019348804\\ 
877& 0.88383663619181803411\\ 
881& 1.75750576898034081030\\ 
883& 0.31716903119407663757\\ 
887& 3.05904637329967833761\\ 
907& 0.31294461576390249253\\ 
911& 3.22665385387042988763\\ 
919& 1.96900000823240599231\\ 
929& 1.69391020650614908554\\ 
937& 2.25537551795571516602\\ 
941& 0.45862903621215601114\\ 
947& 0.51044022041152226105\\ 
953& 1.45295536414148277191\\ 
967& 1.11129489882427547420\\ 
971& 1.51227759095209977084\\ 
977& 1.49820818462619078840\\ 
983& 2.70543359074492325719\\ 
991& 1.69653165274144151606\\ 
997& 0.76277626026681664927\\ 
\phantom{} & \phantom{} \\
\hline
\end{tabular}
}
\caption{\label{TableL}
{\small
Values of $ L(1,\chi_\square)$ for every prime $3\le q \le \boundtable$ with 
$38$-digit precision (we just printed  $20$ digits here);
computed with Pari/GP, v.~2.15.1,
with  the \texttt{lfun}-command.
}
}
\end{table} 
\newpage

\begin{table}[htp]
\scalebox{0.75}{
\begin{tabular}{|r|c|}
\hline
$q$ & $1-c_2(q)/\log q$ \\ \hline
3 & 0.87074071872621259992\\ 
5 & 0.94413465252644938604\\ 
7 & 0.85794172560085927694\\ 
11 & 0.91295899101389699960\\ 
13 & 0.94798100324493175725\\ 
17 & 0.93008660693115370704\\ 
19 & 0.95061813561690229668\\ 
23 & 0.87774529764668030376\\ 
29 & 0.96722478450811694068\\ 
31 & 0.90889864471701991373\\ 
37 & 0.96062352357307734724\\ 
41 & 0.94029273832628963934\\ 
43 & 0.97782499875726672158\\ 
47 & 0.89806478812012098670\\ 
53 & 0.97774365355266156020\\ 
59 & 0.95053773253808441987\\ 
61 & 0.96350960711174903029\\ 
67 & 0.98533494189433124631\\ 
71 & 0.90264833714140339974\\ 
73 & 0.93508619856698228187\\ 
79 & 0.93689666451683397168\\ 
83 & 0.96379164717072082978\\ 
89 & 0.95099544356241617187\\ 
97 & 0.93870159639412537643\\ 
101 & 0.98097919305395174293\\ 
103 & 0.95026245407185662479\\ 
107 & 0.97114985828865582787\\ 
109 & 0.96698365697402200464\\ 
113 & 0.95775902342382907057\\ 
127 & 0.95865616368194021096\\ 
131 & 0.95976323082950450377\\ 
137 & 0.96032886641128995652\\ 
139 & 0.97707497159481733322\\ 
149 & 0.98140111519190970030\\ 
151 & 0.95022029543724605884\\ 
157 & 0.97680962132148529063\\ 
163 & 0.99334329316402472279\\ 
167 & 0.92828921224386521210\\ 
173 & 0.98976456446842987752\\ 
179 & 0.96928960491511367339\\ 
181 & 0.97251506092644502995\\ 
191 & 0.92446210779698297446\\ 
193 & 0.94529931991695971562\\ 
197 & 0.98811334397411123077\\ 
199 & 0.94950144733816672243\\ 
211 & 0.98398396671072490272\\ 
223 & 0.96433881923112980490\\ 
227 & 0.97490737457633045641\\ 
229 & 0.97443605471288035393\\ 
233 & 0.96673597868244165217\\ 
239 & 0.92792073518536729448\\ 
241 & 0.94349611646894449198\\ 
251 & 0.96771964580525059897\\ 
257 & 0.97032902022349032694\\ 
263 & 0.94235143697043577820\\ 
269 & 0.98599209699367395519\\ 
\hline
\end{tabular}
\begin{tabular}{|r|c|}
\hline
$q$ & $1-c_2(q)/\log q$ \\ \hline
271 & 0.95242727966707206534\\ 
277 & 0.97890754448801742017\\ 
281 & 0.96140289769063960078\\ 
283 & 0.98748621922894326509\\ 
293 & 0.99274122047909849711\\ 
307 & 0.98830491300427046313\\ 
311 & 0.92672169478040626367\\ 
313 & 0.95310638797202816364\\ 
317 & 0.98923587662681801345\\ 
331 & 0.98901111738351919822\\ 
337 & 0.95112561283039349772\\ 
347 & 0.98238442972988440478\\ 
349 & 0.97823669495704562692\\ 
353 & 0.97394683751627225384\\ 
359 & 0.93490721444869971096\\ 
367 & 0.96971975879064312552\\ 
373 & 0.98062137611682946929\\ 
379 & 0.99016946005245139097\\ 
383 & 0.94477067621269396065\\ 
389 & 0.98408114281420747217\\ 
397 & 0.98375238754385853422\\ 
401 & 0.96350090660205480320\\ 
409 & 0.94911583765502880144\\ 
419 & 0.97283004226082522827\\ 
421 & 0.97526483902484077226\\ 
431 & 0.93804459164515909435\\ 
433 & 0.95650807631535617704\\ 
439 & 0.95640316150486697055\\ 
443 & 0.98557445348372299620\\ 
449 & 0.96434073755479744957\\ 
457 & 0.95462587889724423098\\ 
461 & 0.98957359647257357008\\ 
463 & 0.98051516124324652035\\ 
467 & 0.98065052483059602863\\ 
479 & 0.93229986446089119777\\ 
487 & 0.98129621017807203356\\ 
491 & 0.97611073843696917771\\ 
499 & 0.99214019247588982038\\ 
503 & 0.94533479957643908790\\ 
509 & 0.98885547089395136569\\ 
521 & 0.96898886485132353691\\ 
523 & 0.98738726867765310271\\ 
541 & 0.97801340786016467275\\ 
547 & 0.99270064115471547721\\ 
557 & 0.99170551183406514282\\ 
563 & 0.97860302943301718019\\ 
569 & 0.96645121794465682260\\ 
571 & 0.98824655917298969209\\ 
577 & 0.96000148436752255790\\ 
587 & 0.98390551443491812454\\ 
593 & 0.97978519297745376467\\ 
599 & 0.94344302545058027601\\ 
601 & 0.95247884344025167162\\ 
607 & 0.97090069047607732686\\ 
613 & 0.98382597591751816930\\ 
617 & 0.97450460801082924052\\ 
\hline
\end{tabular}
\begin{tabular}{|r|c|}
\hline
$q$ & $1-c_2(q)/\log q$ \\ \hline
619 & 0.98898162981970758420\\ 
631 & 0.97178788363623377751\\ 
641 & 0.96606260790537001087\\ 
643 & 0.99358657301462091844\\ 
647 & 0.95107276863659696692\\ 
653 & 0.99008234738545258149\\ 
659 & 0.97693997174342132355\\ 
661 & 0.98092875730547795253\\ 
673 & 0.95793511602798465873\\ 
677 & 0.99486387743516433990\\ 
683 & 0.98981214959834417331\\ 
691 & 0.98990595140998302159\\ 
701 & 0.98727059226367292868\\ 
709 & 0.97818001292129838421\\ 
719 & 0.93935825277518907800\\ 
727 & 0.97479156815475908654\\ 
733 & 0.98792626047388137145\\ 
739 & 0.99042931952965323333\\ 
743 & 0.95997449340728463277\\ 
751 & 0.97165153902162844911\\ 
757 & 0.98236085345455541170\\ 
761 & 0.97377081372221638945\\ 
769 & 0.95530519816588844357\\ 
773 & 0.99422390643906944803\\ 
787 & 0.99089364267401025617\\ 
797 & 0.99325174742430195132\\ 
809 & 0.96895400800129565192\\ 
811 & 0.98754952159246285737\\ 
821 & 0.98917965518610081661\\ 
823 & 0.98417636774549864402\\ 
827 & 0.98773961862416774605\\ 
829 & 0.98088923914296772329\\ 
839 & 0.94285257863898100280\\ 
853 & 0.98892590841388063163\\ 
857 & 0.98207782041867686411\\ 
859 & 0.98810027203440013018\\ 
863 & 0.96443088138037865171\\ 
877 & 0.98612583924837538390\\ 
881 & 0.97244672375952949576\\ 
883 & 0.99500928381600903145\\ 
887 & 0.95192719878892877801\\ 
907 & 0.99511314690650491820\\ 
911 & 0.94967635461707835743\\ 
919 & 0.96936701711950779765\\ 
929 & 0.97383729821145899796\\ 
937 & 0.96524901734749286578\\ 
941 & 0.99294183930919729627\\ 
947 & 0.99212585078030535447\\ 
953 & 0.97771860782924520480\\ 
967 & 0.98295760541235751687\\ 
971 & 0.97683520347803390374\\ 
977 & 0.97718361844355005125\\ 
983 & 0.95870144616958105591\\ 
991 & 0.97416121861100771551\\ 
997 & 0.98844893853186475847\\
\phantom{} & \phantom{} \\
\hline
\end{tabular}
}
\caption{
{\small
Numerical upper bounds in $\beta<1-c_2(q)/\log q$ for every  prime  $3\le q\le \boundtable$ with 
$38$-digit precision (we just printed  $20$ digits here);
computed with Pari/GP, v.~2.15.1.
}
}
\end{table} 

\newpage
\begin{table}[htp]
\scalebox{0.475}{
\begin{tabular}{|c|c|c|}
\hline
$q$  &  $c_{1}(q) $  &  $c_{2}(q)$\\ \hline
3 & 0.55033044351893460168 & 0.14200583483179050854\\ 
5 & 0.26742811116774127267 & 0.08991180821523694362\\ 
7 & 0.61020824229740594980 & 0.27643263791057252254\\ 
11 & 0.39502385106004226257 & 0.20871522398737687604\\ 
13 & 0.25838147218927451454 & 0.13342609230270605368\\ 
17 & 0.35863336447093897006 & 0.19807955817090268479\\ 
19 & 0.24477694706707900762 & 0.14540188635350380661\\ 
23 & 0.62675990410829234931 & 0.38332891209897371249\\ 
29 & 0.18167880710521482528 & 0.11036384645272590779\\ 
31 & 0.49293721594779985979 & 0.31284088835300889574\\ 
37 & 0.22689304729319524630 & 0.14218522406678805450\\ 
41 & 0.34982303782999099771 & 0.22172721913089549226\\ 
43 & 0.12737646854812820283 & 0.08340461723967329728\\ 
47 & 0.59510495844651555500 & 0.39246561154912285570\\ 
53 & 0.13601643326508233304 & 0.08836419232508230269\\ 
59 & 0.30091730507159657222 & 0.20168424763644029217\\ 
61 & 0.22825078784983049459 & 0.15000740241772657957\\ 
67 & 0.09128054381737623086 & 0.06166206157984507066\\ 
71 & 0.61226018665265819509 & 0.41497897426385014935\\ 
73 & 0.41828538606457288375 & 0.27851003221912283015\\ 
79 & 0.40446390522732929947 & 0.27572673371042579858\\ 
83 & 0.23411203757965752099 & 0.15999893982344563455\\ 
89 & 0.32625529930182564596 & 0.21996363430833003744\\ 
97 & 0.41390490735809004215 & 0.28042247994052598910\\ 
101 & 0.12928560930103269329 & 0.08778331638358404990\\ 
103 & 0.33394651871870195182 & 0.23052004591656898579\\ 
107 & 0.19498417508113876678 & 0.13481177406708121097\\ 
109 & 0.22722075211365497770 & 0.15489115093387282650\\ 
113 & 0.29241343466467224273 & 0.19968947811670376211\\ 
127 & 0.28773792602782610679 & 0.20027727799660294114\\ 
131 & 0.28150883583268095963 & 0.19616218935426226200\\ 
137 & 0.28329650751350412613 & 0.19518122056243286339\\ 
139 & 0.16200292957965930315 & 0.11312295508165454357\\ 
149 & 0.13461292190557461195 & 0.09306782093014857491\\ 
151 & 0.35669014859897260709 & 0.24975870798530923793\\ 
157 & 0.16924788534457484666 & 0.11725625491767883858\\ 
163 & 0.04830792988513398745 & 0.03390760178246081580\\ 
167 & 0.52249791984297772527 & 0.36701536801942904653\\ 
173 & 0.07585653441526830736 & 0.05274618389086421808\\ 
179 & 0.22633140124075549858 & 0.15930666755510101877\\ 
181 & 0.20515012815743757396 & 0.14288037417839760684\\ 
191 & 0.56263819926509207734 & 0.39674566402855702698\\ 
193 & 0.41241911293129269104 & 0.28787273239944098485\\ 
197 & 0.08990773646468201719 & 0.06279962543820142560\\ 
199 & 0.37865085304891402689 & 0.26730423244648594619\\ 
211 & 0.12123424357604100845 & 0.08571553802523036784\\ 
223 & 0.27234870527388714480 & 0.19282612999037137270\\ 
227 & 0.19218133559014944790 & 0.13612623873079037965\\ 
229 & 0.19792483225116621922 & 0.13890737200427872534\\ 
233 & 0.25822880601066285047 & 0.18132345932223975371\\ 
239 & 0.55659843343745476824 & 0.39473946660735538723\\ 
241 & 0.44091994894713225232 & 0.30991232712144637333\\ 
251 & 0.25121361182915890473 & 0.17836357796159308676\\ 
257 & 0.23381999795655663869 & 0.16464652429323954028\\ 
263 & 0.45195192190267889927 & 0.32122667293444028461\\ 
269 & 0.11115733759606550144 & 0.07837017435385113565\\ 
\hline
\end{tabular}
}
\scalebox{0.475}{
\begin{tabular}{|c|c|c|}
\hline
$q$  &  $c_{1}(q) $  &  $c_{2}(q)$\\ \hline
271 & 0.37471889822313975596 & 0.26650803193754201335\\ 
277 & 0.16812095224667437033 & 0.11862433904786765068\\ 
281 & 0.30831279144623607124 & 0.21762415202873450876\\ 
283 & 0.09923836100529680835 & 0.07064588483174984608\\ 
293 & 0.05835005678626515758 & 0.04123112060951874748\\ 
307 & 0.09392610270553316042 & 0.06697598261933001858\\ 
311 & 0.58969447054985015578 & 0.42060229691593070488\\ 
313 & 0.38071266883740373676 & 0.26946022305108387068\\ 
317 & 0.08755628574697731273 & 0.06198952918805163974\\ 
331 & 0.08928336406008632715 & 0.06375879775394483775\\ 
337 & 0.40119076894989782797 & 0.28445298649725928204\\ 
347 & 0.14416153690968422288 & 0.10303919169388233892\\ 
349 & 0.17957640093915020960 & 0.12742571629132863783\\ 
353 & 0.21533648778321692919 & 0.15284004549288177497\\ 
359 & 0.53546809878953853528 & 0.38296184256303134496\\ 
367 & 0.24992680428189637024 & 0.17881578118762586811\\ 
373 & 0.16147759244702068455 & 0.11475204098897685117\\ 
379 & 0.08153522079748106318 & 0.05836918685407896260\\ 
383 & 0.45880262340283667315 & 0.32850595031568327195\\ 
389 & 0.13346901260038062705 & 0.09493336788720439026\\ 
397 & 0.13663209108292684910 & 0.09722467765085002936\\ 
401 & 0.30738441205713295880 & 0.21877415795894336801\\ 
409 & 0.42976883702961221355 & 0.30600285829649572809\\ 
419 & 0.22877130332366804116 & 0.16404869772887678440\\ 
421 & 0.20979552591959329948 & 0.14946549585491626113\\ 
431 & 0.52386621109724941649 & 0.37582820384698174496\\ 
433 & 0.37039307941713669135 & 0.26402806197577231845\\ 
439 & 0.36964430375980270234 & 0.26526493823556995421\\ 
443 & 0.12247465530556111599 & 0.08790307416796546949\\ 
449 & 0.30528681149248847923 & 0.21777193191285979942\\ 
457 & 0.38944967105004687527 & 0.27790212589447049403\\ 
461 & 0.08960317341694861572 & 0.06394928299061061816\\ 
463 & 0.16651365196643178377 & 0.11959262198183381426\\ 
467 & 0.16556695088817405610 & 0.11892824535424547236\\ 
479 & 0.58145656097212776090 & 0.41782496695151772864\\ 
487 & 0.16103289122396401185 & 0.11574399152071413545\\ 
491 & 0.20592501105759230932 & 0.14802847452958968898\\ 
499 & 0.06791208062003553977 & 0.04882988813572055338\\ 
503 & 0.47288287326840516245 & 0.34004980840133510986\\ 
509 & 0.09714492157105550040 & 0.06945769832238238121\\ 
521 & 0.27121802009494090607 & 0.19399791000115208837\\ 
523 & 0.10972958501609965330 & 0.07895041919659370235\\ 
541 & 0.19332299637930880256 & 0.13837084284897416421\\ 
547 & 0.06391905447241206028 & 0.04601843413060230529\\ 
557 & 0.07323316900168368336 & 0.05244244256092793828\\ 
563 & 0.18815229551239650389 & 0.13551299779577685311\\ 
569 & 0.29709999451045999039 & 0.21282946206966013679\\ 
571 & 0.10356349389164450671 & 0.07460366348169564142\\ 
577 & 0.35491533119298588274 & 0.25430425328575003191\\ 
587 & 0.14237934219934102363 & 0.10260274493976198215\\ 
593 & 0.18006195764727533672 & 0.12907547257718210897\\ 
599 & 0.50178555632224349672 & 0.36169664744149677984\\ 
601 & 0.42408857717593404003 & 0.30406863164646020543\\ 
607 & 0.25866671023500976377 & 0.18648376288401031238\\ 
613 & 0.14474089554896833741 & 0.10381078904398923563\\ 
617 & 0.22836539873653918938 & 0.16380455424354878123\\ 
\hline
\end{tabular}
}
\scalebox{0.475}{
\begin{tabular}{|c|c|c|}
\hline
$q$  &  $c_{1}(q) $  &  $c_{2}(q)$\\ \hline
619 & 0.09821810981007854740 & 0.07082724345232840177\\ 
631 & 0.25217397926401117719 & 0.18189214322680485918\\ 
641 & 0.30560600168191523306 & 0.21933836479866029927\\ 
643 & 0.05748045804689990168 & 0.04147014706579212793\\ 
647 & 0.43889822330321299109 & 0.31667398461510986618\\ 
653 & 0.08953939265010298985 & 0.06428203036255914768\\ 
659 & 0.20739903882064491981 & 0.14967626811126803621\\ 
661 & 0.17247239110369827858 & 0.12384395547829596586\\ 
673 & 0.38136874412082350586 & 0.27391581174681826550\\ 
677 & 0.04660340526027440933 & 0.03347555849498589539\\ 
683 & 0.09209361584897920724 & 0.06649095327648296169\\ 
691 & 0.09139588640955929171 & 0.06599630106943590398\\ 
701 & 0.11606038625471207803 & 0.08340954458956226259\\ 
709 & 0.19925564553833806973 & 0.14322324277539521284\\ 
719 & 0.55215677551823390162 & 0.39889300573477128647\\ 
727 & 0.22988538302098122046 & 0.16609650404220700767\\ 
733 & 0.11076160532740114243 & 0.07965221882043248666\\ 
739 & 0.08747928954797630044 & 0.06321719581284105574\\ 
743 & 0.36612347141509609398 & 0.26459645812028807150\\ 
751 & 0.25969941304500925429 & 0.18770665974100441980\\ 
757 & 0.16253453608508091490 & 0.11693630993036847128\\ 
761 & 0.24186135198530417594 & 0.17402103422815112884\\ 
769 & 0.41272464112715707084 & 0.29700102405169899696\\ 
773 & 0.05337597163560077649 & 0.03841263399093716531\\ 
787 & 0.08396959410194464642 & 0.06072326916134757729\\ 
797 & 0.06262040466503263606 & 0.04508409479401042827\\ 
809 & 0.28867828366497874651 & 0.20787771960393117173\\ 
811 & 0.11528557294782929567 & 0.08339664177565978132\\ 
821 & 0.10081365428925087305 & 0.07261017392591420386\\ 
823 & 0.14681777932965260644 & 0.10622335026008855922\\ 
827 & 0.11383289589061803950 & 0.08236284756934693961\\ 
829 & 0.17829052679829539376 & 0.12842852029139896932\\ 
839 & 0.53164915993572322425 & 0.38472848193364481224\\ 
853 & 0.10371429843127831870 & 0.07473638132162035954\\ 
857 & 0.16795630454669542980 & 0.12103632716842654726\\ 
859 & 0.11106482700547295177 & 0.08039181237006652880\\ 
863 & 0.33219296634601661470 & 0.24046199206497726227\\ 
877 & 0.13042658071284927743 & 0.09401834734667991368\\ 
881 & 0.25917871015551250574 & 0.18684035396996865892\\ 
883 & 0.04675716139421373779 & 0.03385365085992796650\\ 
887 & 0.45066532622240747506 & 0.32631072248613567032\\ 
907 & 0.04595272684180220806 & 0.03328016569948733107\\ 
911 & 0.47349527364011849896 & 0.34293264020841757477\\ 
919 & 0.28857063487066901495 & 0.20901760697482178406\\ 
929 & 0.24786117272113414772 & 0.17879874892404850467\\ 
937 & 0.32960396162277473481 & 0.23778996803766798800\\ 
941 & 0.06698303562075857136 & 0.04832682491998277619\\ 
947 & 0.07448094902483237890 & 0.05396389970694517977\\ 
953 & 0.21181296392709675170 & 0.15284176980867364772\\ 
967 & 0.16166174130108460363 & 0.11715280323329332387\\ 
971 & 0.21986126973234104798 & 0.15933503304962222559\\ 
977 & 0.21762090048960129350 & 0.15707907427329480398\\ 
983 & 0.39262618784087143443 & 0.28457219167277347527\\ 
991 & 0.24591996730680940513 & 0.17825437671761626813\\ 
997 & 0.11047122273915650913 & 0.07975720056598333162\\
 \phantom{} & \phantom{} & \phantom{}  \\
\hline
\end{tabular}
}
\caption{
{\small
Values of $c_{1}(q)$ and $c_{2}(q)$ for every prime  $3\le q\le \boundtable$ with 
$38$-digit precision  (we just printed  $20$ digits here);
computed with Pari/GP, v.~2.15.1,
with  the \texttt{lfun}-command. 
}
}
\end{table} 

\newpage
\begin{table}[htp]
\scalebox{0.475}{
\begin{tabular}{|c|c|c|}
\hline
$q$  &  $c_{3}(q) $  &  $c_{4}(q)$\\ \hline
3 & 1.60531281426088943415 & 2.27009444749574814599\\ 
5 & 1.44948863673327054404 & 1.52484940046373303805\\ 
7 & 1.40549371558030080078 & 0.80194548881209684904\\ 
11 & 1.37578983698684533502 & 0.51685528737616402827\\ 
13 & 1.36949704631674633504 & 0.56701677740432957764\\ 
17 & 1.36252236893620521850 & 0.44802975356757556720\\ 
19 & 1.36046423241498071674 & 0.32298674064855080025\\ 
23 & 1.35776141486577218049 & 0.27728327997964533856\\ 
29 & 1.35554642561250613691 & 0.29063403062083574784\\ 
31 & 1.35507464784761415630 & 0.22060562381678128245\\ 
37 & 1.35410022350787716690 & 0.24165664069777158523\\ 
41 & 1.35368825401096606619 & 0.22402980514835986080\\ 
43 & 1.35352908070133037244 & 0.17368215488406234444\\ 
47 & 1.35327945785782956839 & 0.16304437010137342127\\ 
53 & 1.35302803290054631125 & 0.18624290576986660697\\ 
59 & 1.35287356312845265771 & 0.13914829082434390134\\ 
61 & 1.35283715648703579283 & 0.16875967254290138130\\ 
67 & 1.35276058597980783655 & 0.12757499604056836693\\ 
71 & 1.35273045601505508280 & 0.12267004527244762442\\ 
73 & 1.35272021503884114050 & 0.14914807577906107908\\ 
79 & 1.35270458240959514828 & 0.11419672158984247292\\ 
83 & 1.35270405753963782554 & 0.11050587271361023816\\ 
89 & 1.35271409234213129784 & 0.13050971566471038628\\ 
97 & 1.35274181188616268654 & 0.12326292013772485509\\ 
101 & 1.35275999626165793191 & 0.12002110414809822468\\ 
103 & 1.35276990437852579440 & 0.09589594761234884346\\ 
107 & 1.35279105760962685026 & 0.09355275345925109844\\ 
109 & 1.35280220227926316358 & 0.11416831697694400180\\ 
113 & 1.35282541291606598287 & 0.11151530960201366974\\ 
127 & 1.35291296473365380482 & 0.08378484195766069264\\ 
131 & 1.35293899675470203652 & 0.08214304815603021596\\ 
137 & 1.35297838607412964315 & 0.09847532805450864286\\ 
139 & 1.35299155829969711774 & 0.07910442489595189351\\  
149 & 1.35305731750500008284 & 0.09333866082791299003\\ 
151 & 1.35307039776105200217 & 0.07506859077004280622\\ 
157 & 1.35310939233420632015 & 0.09029237299931346840\\ 
163 & 1.35314792909217273090 & 0.07154527672817944797\\ 
167 & 1.35317332141301258971 & 0.07046711812116448849\\ 
173 & 1.35321091223368340488 & 0.08493169486413290901\\ 
179 & 1.35324786293024902366 & 0.06747987379471857938\\ 
181 & 1.35326003133792057582 & 0.08255737421626130887\\ 
191 & 1.35331972084050707793 & 0.06481348205845241569\\ 
193 & 1.35333142439847546151 & 0.07931247172099445145\\ 
197 & 1.35335459463309306442 & 0.07830579883782960092\\ 
199 & 1.35336606106221824264 & 0.06318783927574056066\\ 
211 & 1.35343320063344513445 & 0.06094564326246216244\\ 
223 & 1.35349753210383306421 & 0.05890806503949326757\\ 
227 & 1.35351836565596155042 & 0.05826923225318402691\\ 
229 & 1.35352867011805046692 & 0.07134050273671778391\\ 
233 & 1.35354905701723639571 & 0.07058440577838804781\\ 
239 & 1.35357909051582260329 & 0.05646089920504149769\\ 
241 & 1.35358895812024356838 & 0.06913582679722596790\\ 
251 & 1.35363724585979693559 & 0.05479834788710899596\\ 
257 & 1.35366539993686259322 & 0.06646781550907920780\\ 
263 & 1.35369296016889736324 & 0.05326343536964672312\\ 
269 & 1.35371994301608376598 & 0.06464282716452377542\\ 
\hline
\end{tabular}
}
\scalebox{0.475}{
\begin{tabular}{|c|c|c|}
\hline
$q$  &  $c_{3}(q) $  &  $c_{4}(q)$\\ \hline
271 & 1.35372881182136547397 & 0.05230347733289799644\\ 
277 & 1.35375505021834899060 & 0.06350007293840636272\\ 
281 & 1.35377224174122466851 & 0.06294914957533875233\\ 
283 & 1.35378074904363092690 & 0.05094878674805416513\\ 
293 & 1.35382242403322780813 & 0.06137211649047258609\\ 
307 & 1.35387845858578453571 & 0.04850638195293471844\\ 
311 & 1.35389399608837588049 & 0.04812989886842535723\\ 
313 & 1.35390168872018058181 & 0.05896980870015960658\\ 
317 & 1.35391692427449962776 & 0.05851977098360236817\\ 
331 & 1.35396872821695014793 & 0.04636137845546200668\\ 
337 & 1.35399023422861369381 & 0.05640370780779139758\\ 
347 & 1.35402519507511206504 & 0.04506901888157027469\\ 
349 & 1.35403205874034062285 & 0.05523136284868475601\\ 
353 & 1.35404566089891669524 & 0.05485530539427773666\\ 
359 & 1.35406575725457480815 & 0.04416257590215822297\\ 
367 & 1.35409199598600849378 & 0.04358556145315234400\\ 
373 & 1.35411127094265687404 & 0.05307583488435248998\\ 
379 & 1.35413021042820747073 & 0.04275786000725053118\\ 
383 & 1.35414265504853999390 & 0.04249147885966666124\\ 
389 & 1.35416105619981949696 & 0.05176201985666864596\\ 
397 & 1.35418510867177370227 & 0.05113804985383631445\\ 
401 & 1.35419693385708307854 & 0.05083378332121723922\\ 
409 & 1.35422019428835176059 & 0.05024001050227069348\\ 
419 & 1.35424856267586232283 & 0.04028431988171994174\\ 
421 & 1.35425414487149731846 & 0.04938436527159565703\\ 
431 & 1.35428161319598569452 & 0.03961646506144404360\\ 
433 & 1.35428702014259397731 & 0.04856795158054284160\\ 
439 & 1.35430307207152657818 & 0.03918792684646118666\\ 
443 & 1.35431363503286562166 & 0.03897842516898990112\\ 
449 & 1.35432927673113316805 & 0.04753554951583818809\\ 
457 & 1.35434976349460983667 & 0.04704171487765358980\\ 
461 & 1.35435985279916674685 & 0.04680008561804179954\\ 
463 & 1.35436485968688543723 & 0.03797564767087865659\\ 
467 & 1.35437479898402282743 & 0.03778355999124213694\\ 
479 & 1.35440403571820347686 & 0.03722309756317605641\\ 
487 & 1.35442305829226240461 & 0.03686204522125124301\\ 
491 & 1.35443243353832925209 & 0.03668513159221982500\\ 
499 & 1.35445091917610656679 & 0.03633827188812157807\\ 
503 & 1.35446003236494670358 & 0.03616822733487120387\\ 
509 & 1.35447354392782473163 & 0.04414639190945388174\\ 
521 & 1.35450001358879804956 & 0.04354607924452285284\\ 
523 & 1.35450435529392795671 & 0.03535001321630170097\\ 
541 & 1.35454256786249600342 & 0.04259423060531816795\\ 
547 & 1.35455497140447127621 & 0.03443306523769305016\\ 
557 & 1.35457528831888197854 & 0.04187318786719491501\\ 
563 & 1.35458727076482901904 & 0.03385736977031018334\\ 
569 & 1.35459910155248485596 & 0.04135422004740572241\\ 
571 & 1.35460301201156021611 & 0.03357940502288334637\\ 
577 & 1.35461464563153277280 & 0.04101805263161343844\\ 
587 & 1.35463371601799856365 & 0.03304204531696241530\\ 
593 & 1.35464497175843319268 & 0.04036799276590968032\\ 
599 & 1.35465609117902202183 & 0.03265440747138537032\\ 
601 & 1.35465976784663745689 & 0.04005357328162699404\\ 
607 & 1.35467070985905401855 & 0.03240292183238713396\\ 
613 & 1.35468152197615147700 & 0.03959451505151431097\\ 
617 & 1.35468865918922608328 & 0.03944473222429774966\\ 
\hline
\end{tabular}
}
\scalebox{0.475}{
\begin{tabular}{|c|c|c|}
\hline
$q$  &  $c_{3}(q) $  &  $c_{4}(q)$\\ \hline
619 & 1.35469220682238232404 & 0.03203562589619692571\\ 
631 & 1.35471320483555636308 & 0.03167971346422158052\\ 
641 & 1.35473033619057541217 & 0.03857813477593126179\\ 
643 & 1.35473372347150569529 & 0.03133461999753278949\\ 
647 & 1.35474045971790388657 & 0.03122190102561039225\\ 
653 & 1.35475046953239467792 & 0.03816434888713844033\\ 
659 & 1.35476036781470475868 & 0.03089041996727199247\\ 
661 & 1.35476364282570692151 & 0.03789529186024628706\\ 
673 & 1.35478304188845849460 & 0.03750147680637506693\\ 
677 & 1.35478941441495095188 & 0.03737272881724802515\\ 
683 & 1.35479888724152173116 & 0.03025594068424354252\\ 
691 & 1.35481136016874090170 & 0.03005241141445773599\\ 
701 & 1.35482670433578210082 & 0.03662539896746132327\\ 
709 & 1.35483878722270831350 & 0.03638544200950989301\\ 
719 & 1.35485365694129056644 & 0.02936909848443259670\\ 
727 & 1.35486537037338342442 & 0.02918171927857417683\\ 
733 & 1.35487405156457705708 & 0.03569117019153599530\\ 
739 & 1.35488264543316543449 & 0.02890681917650937155\\ 
743 & 1.35488832691210862580 & 0.02881678394259856258\\ 
751 & 1.35489957719472852988 & 0.02863904313626900627\\ 
757 & 1.35490791799401157279 & 0.03503277862848840545\\ 
761 & 1.35491343315962672207 & 0.03492632418953755422\\ 
769 & 1.35492435643738726104 & 0.03471610842280814580\\ 
773 & 1.35492976534440364969 & 0.03461232250601469780\\ 
787 & 1.35494842626122993014 & 0.02787557616244300825\\ 
797 & 1.35496150462001991354 & 0.03400737567237013774\\ 
809 & 1.35497693157006976121 & 0.03371582595582892482\\ 
811 & 1.35497947500667002529 & 0.02739720699393118373\\ 
821 & 1.35499207589280056863 & 0.03343118270059318275\\ 
823 & 1.35499457309647581098 & 0.02716649572848899616\\ 
827 & 1.35499954485148919285 & 0.02709079381454957239\\ 
829 & 1.35500201948179621086 & 0.03324512681998789564\\ 
839 & 1.35501428171472255976 & 0.02686718768846612490\\ 
853 & 1.35503114514538460618 & 0.03270393912307176912\\ 
857 & 1.35503589973501215812 & 0.03261612570176427643\\ 
859 & 1.35503826661981379802 & 0.02650574552671867219\\ 
863 & 1.35504297974445551695 & 0.02643508044000786574\\ 
877 & 1.35505926250729193789 & 0.03218678460645850450\\ 
881 & 1.35506385492579003853 & 0.03210280713893170195\\ 
883 & 1.35506614133270666050 & 0.02608951563543007725\\ 
887 & 1.35507069470268753185 & 0.02602191172110351678\\ 
907 & 1.35509308082340970970 & 0.02569111528743156599\\ 
911 & 1.35509748361184231270 & 0.02562636160846504475\\ 
919 & 1.35510621664801272634 & 0.02549821685997481607\\ 
929 & 1.35511699928640912613 & 0.03114087003695233544\\ 
937 & 1.35512552064462352496 & 0.03098830216628082727\\ 
941 & 1.35512974700391999364 & 0.03091280248419242123\\ 
947 & 1.35513604424031454818 & 0.02506348598834765208\\ 
953 & 1.35514229133818874919 & 0.03068937098973946392\\ 
967 & 1.35515667676797795161 & 0.02476541505903528806\\ 
971 & 1.35516073863466659129 & 0.02470698740314042627\\ 
977 & 1.35516679194687294594 & 0.03025581444672276531\\ 
983 & 1.35517279844033022223 & 0.02453400120490038280\\ 
991 & 1.35518073535005806431 & 0.02442055038098664765\\ 
997 & 1.35518663501537123042 & 0.02990739982056048090\\
 \phantom{} & \phantom{} & \phantom{}  \\
\hline
\end{tabular}
}
\caption{\label{Tablec3c4}
{\small
Values of $c_{3}(q)$ and $c_{4}(q)$ for every prime $3\le q \le \boundtable$  with 
$38$-digits precision (we just printed  $20$ digits here);
computed with Pari/GP, v.~2.15.1. 
}
}
\end{table}

\newpage  
\restoregeometry   
\begin{figure}[ht] 
\includegraphics[scale=0.45,angle=0]{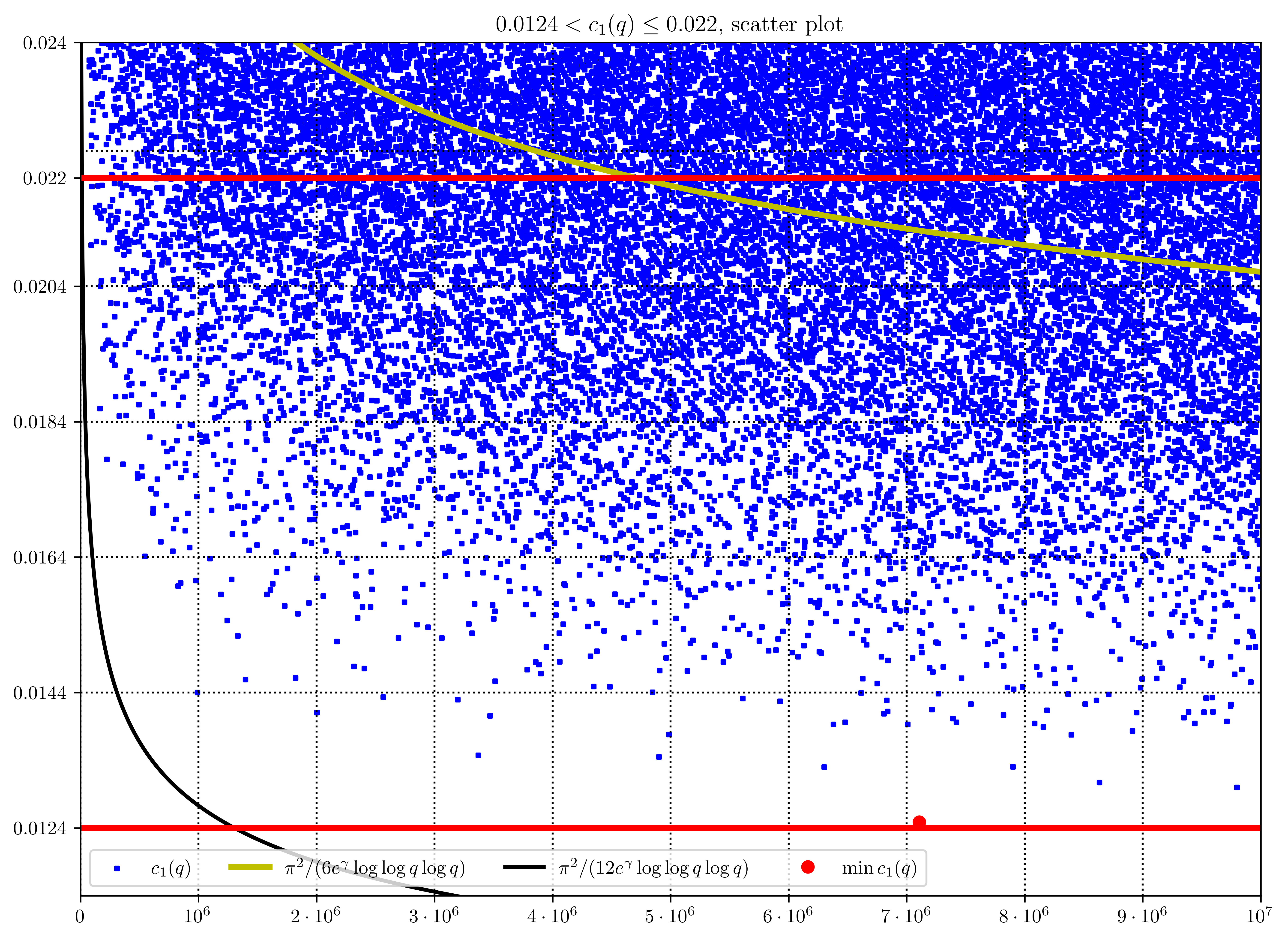}  
\caption{\small{The values of $0.0124<c_1(q)\le 0.022$, $q$ prime, $3\le q\le \bound$, in Theorem \ref{Thm-beta}.
$ \min c_{1}(q) = \constonemin\dotsc$ attained at 
$q=\constoneminattained$.
The red straight lines correspond to $0.0124$ and $0.022$.
The yellow one corresponds to the second Joshi bound, see \eqref{Joshi-second};
the black one corresponds to the second Littlewood bound, see \eqref{Littlewood-bounds}.
 }}
\end{figure}

\begin{figure}[ht] 
\includegraphics[scale=0.5,angle=0]{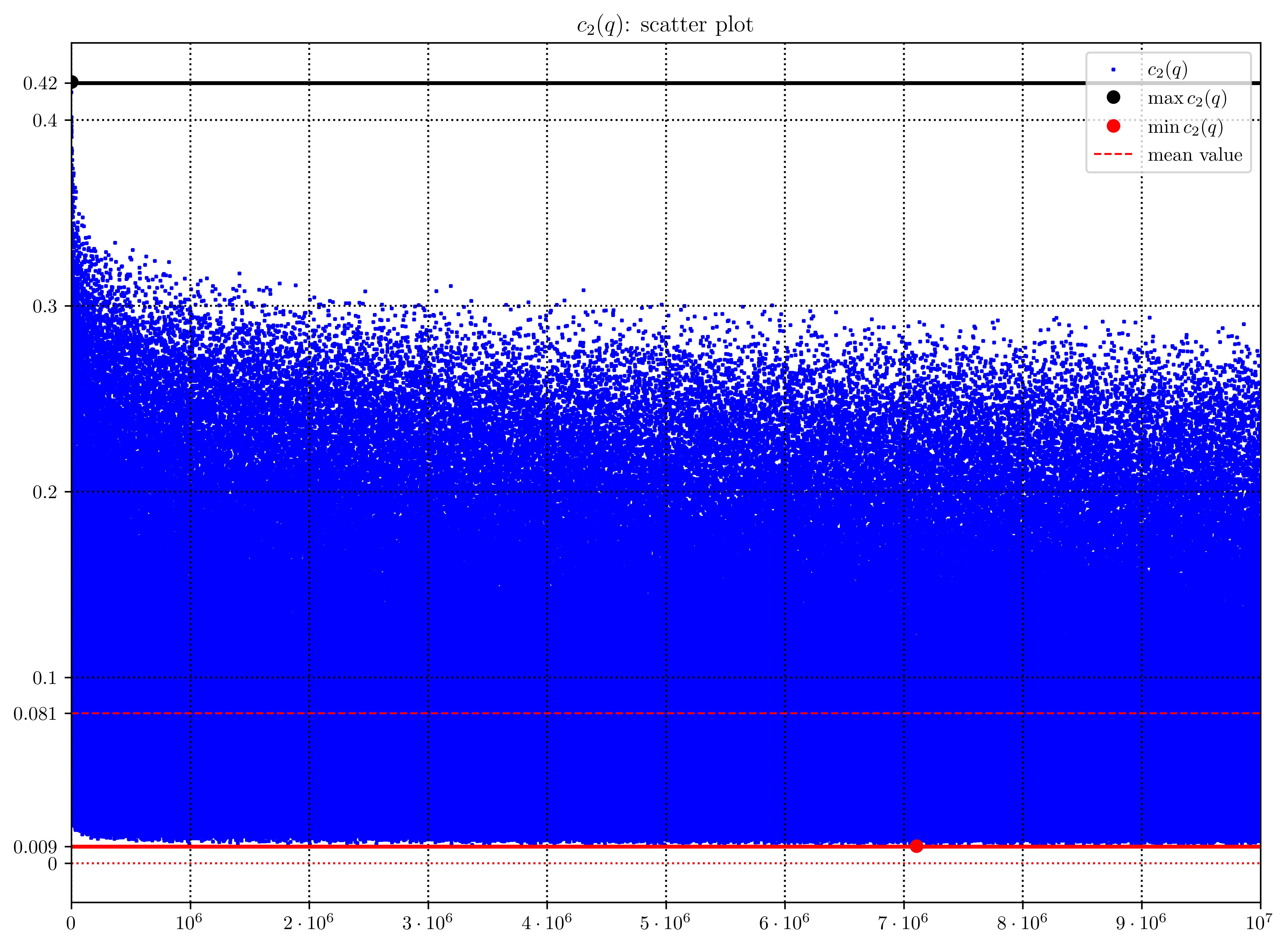}  
\caption{\small{The values of $c_2(q)$, $q$ prime, $3\le q\le \bound$, in Theorem \ref{Thm-beta}.
$ \min c_{2}(q) = \consttwomin\dotsc$ attained at 
$q=\consttwominattained$; 
$ \max c_{2}(q) =  \consttwomax\dotsc$ attained at 
$q=\consttwomaxattained$.
The black line corresponds to $0.42$;  the red one  to $0.009$.
The red dashed line corresponds to the mean value.
 }}
 \end{figure}

\end{document}